\documentclass[a4paper]{amsart}

\usepackage{txfonts}
\usepackage{hyperref}
\usepackage{amssymb}
\usepackage{amsthm}
\usepackage{amscd}
\usepackage{amsmath}
\usepackage[all]{xy}
\usepackage{graphicx}
\theoremstyle{plain}
\newtheorem{theorem}{Theorem}
\newtheorem*{theorem*}{Theorem}

\newtheorem*{maintheorem*}{Main Theorem}
\newtheorem{proposition}[theorem]{Proposition}%[chapter]
\newtheorem{lemma}[theorem]{Lemma}%[chapter]
\newtheorem{claim}[theorem]{Claim}
\newtheorem*{conjecture*}{Conjecture}%[chapter]
\newtheorem{fact}[theorem]{Fact}
\theoremstyle{definition}

\newtheorem*{definition*}{Definition}%[chapter]
\newtheorem{example}{Example}%[chapter]
\newtheorem*{example*}{Example}
\newtheorem{notation}{Notation}%[chapter]
\newtheorem*{notation*}{Notation}
\newtheorem*{notation-conv*}{Notation and convention}
\newtheorem*{convention*}{Convention}
\theoremstyle{remark}
\newtheorem{remark}{Remark}%[chapter]

\newcommand{\Gcenter}[2]{
  \dimen0=\ht\strutbox%
  \advance\dimen0\dp\strutbox%
  \multiply\dimen0 by#1%
  \divide\dimen0 by2%
  \advance\dimen0 by-.5\normalbaselineskip
  \raisebox{-\dimen0}[0pt][0pt]{#2}
}
\newcommand{\Z}{{\mathbb Z}}
\newcommand{\ZZ}{{\mathbb Z}}
\newcommand{\C}{{\mathbb C}}
\newcommand{\CC}{{\mathbb C}}

\newcommand{\IR}{{\mathbb R}}

\newcommand{\cusp}{{\mathfrak{c}}}

\newcommand{\I}{\mathbf{1}}

\newcommand{\bord}{\partial}

\newcommand{\SL}{{\mathrm{SL}_2(\C)}}
\newcommand{\PSL}{{\mathrm{PSL}_2(\C)}}

\newcommand{\T}{\mathcal{T}}
\newcommand{\sll}{\mathfrak{sl}_2(\CC)}
\newcommand{\tsll}{\widetilde{\mathfrak{sl}}_2(\CC)}

\newcommand{\trace}{{\rm tr}\,}
\newfont{\bg}{cmr10 scaled\magstep4}

\newcommand{\tr}{\mathop{\mathrm{tr}}\nolimits}
\newcommand{\im}{\mathop{\mathrm{im}}\nolimits}

\begin{document}

%%%%%%%%%%%%%%%%%%%%%%%%%
% Subject classification 
%%%%%%%%%%%%%%%%%%%%%%%%%
%
%

%%%%%%%%%%%%%%%%%%%%%%%
%  title 
%%%%%%%%%%%%%%%%%%%%%%%

\title{Non--abelian Reidemeister torsion for twist knots}

%%%%%%%%%%%%%%%%%%%%%%%
%  author names and addresses
%%%%%%%%%%%%%%%%%%%%%%%

\author{J\'er\^ome Dubois \and Vu Huynh \and Yoshikazu Yamaguchi}

\address{Centre de Recerca Matematica, Apartat 50, E-08193 Bellaterra, Spain \\ \& Institut de Math\'ematiques de Jussieu, 
Universit\'e Paris 7 (Denis Diderot), 
Case 7012, 
2, place Jussieu, 
75251 Paris Cedex 05
France, dubois@math.jussieu.fr}

\address{Faculty of Mathematics and Informatics, University of Natural Sciences, Vietnam National University, 227 Nguyen Van Cu, District 5, Ho Chi Minh City, Vietnam}
\email{hqvu@hcmuns.edu.vn}

\address{Graduate School of Mathematical Sciences, 
University of Tokyo, 3-8-1 Komaba Meguro, 
Tokyo 153-8914, Japan}
\email{shouji@ms.u-tokyo.ac.jp}

\date{\today}

%%%%%%%%%%%%%%%%%%%%%%%%%
%  abstract
%%%%%%%%%%%%%%%%%%%%%%%%%
\begin{abstract}
    This paper gives an explicit formula for the $\mathrm{SL}_{\;2}(\mathbb{C})-$non--abelian Reidemeister torsion as defined in~\cite{JDFibre} in the case of twist knots. For hyperbolic twist knots, we also prove that the non--abelian Reidemeister torsion at the holonomy representation can be expressed as a rational function evaluated at the cusp shape of the knot. 
\end{abstract}

\maketitle

%\tableofcontents

%%%%%%%%%%%%%%%%%%%%%%%%%
% end Topmatter
%%%%%%%%%%%%%%%%%%%%%%%%%

%%%%%%%%%%%%%%%%%%%%%
% body of paper
%%%%%%%%%%%%%%%%%%%%%

%%%%%%%%%%%%%%%%%%%%%%%%%%%%%%%%%%%
\section{Introduction}
%%%%%%%%%%%%%%%%%%%%%%%%%%%%%%%%%%%
    Twist knots form a family of special two--bridge knots which include the trefoil knot and the figure eight knot. The knot group of a two--bridge knot has a particularly nice presentation  with only two generators and a single relation. One could find our interest in this family of knots in the following facts: first, twist knots except the trefoil knot are hyperbolic; and second, twist knots are not fibered except the trefoil knot and the figure eight knot (see Remark~\ref{remarktwist} of the present paper for details).

The non--abelian Reidemeister torsion associated to a representation of a knot group to a general linear group over a field has been studied since the early 1990's. {It was} initially considered as a twisted Alexander polynomial by Lin  \cite{Lin01},  Wada \cite{Wada94}, and later interpreted as a form of Reidemeister torsion by Kitano \cite{Kitano}, Kirk-Livingston \cite{KL} and Goda--Kitano and Morifuji~\cite{GKM}. This invariant in many cases is stronger than the classical ones,  for example  it detects the unknot \cite{SW06}, and decides fiberness for knots of genus one  \cite{FV07}. Abelian Reidemeister torsions are now well--known, see e.g. Turaev's monograph~\cite{Turaev:2002}; but unfortunately in the case of non--abelian representations concrete computations of such torsions are still very few. 
In \cite{Porti:1997}, Porti began the study of the non--abelian Reidemeister torsion (consider as a functional on the non--abelian part of the character variety) with the adjoint representation associated to an irreducible representation of the fundamental group of a hyperbolic three--dimensional manifold to $\SL$. In~\cite{JDFourier}, the first author introduced a sign--refined version of this torsion. 
%the notion of the a sign-refined non--abelian Reidemeister torsion in the adjoint representation associated to an irreducible representation of a knot group to $\SL$. 
In the present paper, we call this (sign--refined) torsion  the $\SL-$non--abelian Reidemeister torsion. One can observe that this torsion has connections with hyperbolic structures, the theory of character variety, and the theory of Chern-Simons invariant, see e.g. \cite{DK05}.

In~\cite[Main Theorem]{JDFibre},  one can find an \lq\lq explicit\rq\rq\, formula which gives the value of {the non--abelian Reidemeister} torsion for fibered knots in terms of the map induced by the monodromy of the knot at the level of the character variety of the knot exterior. In particular, a practical formula of the non--abelian Reidemeister torsion for torus knots is presented in~\cite[Section 6.2]{JDFibre}.
One can also find an explicit formula for the non--abelian Reidemeister torsion for the figure eight knot in~\cite[Section 7]{JDFibre}. More recently, the last author found~\cite[Theorem 3.1.2]{YY1} an interpretation of the non--abelian Reidemeister torsion in terms of a sort of twisted Alexander polynomial (called in this paper, the non--abelian Reidemeister torsion polynomial) and gave an explicit formula of the non--abelian torsion for the twist knot $5_2$. 
    
    In the present paper we give an explicit formula of the non--abelian Reidemeister torsion for all twist knots. Since twist knots are particular two--bridge knots, this paper is a first step in the understanding of the non--abelian Reidemeister torsion for two--bridge knots.

\section*{Organization} 

The outline of the paper is as follows. We recall some properties of twist knots in Section~\ref{section:twist_knots}.
In Section~\ref{section:review_character_variety}, we give a review of the character variety of twist knots including remarks on parabolic representations and Riley's method for describing the non--abelian part of the character variety. Next we give an explicit formula for
the cusp shape of hyperbolic twist knots.
In Section~\ref{section:review_torsion}, we recall the definition of 
the non--abelian Reidemeister torsion for a knot and an algebraic
description of this invariant.
We give formulas for the non--abelian Reidemeister torsion for twist knots 
in Section~\ref{section:results}.
In particular, we show in Subsection~\ref{subsection:results_holonomy} that the non--abelian Reidemeister torsion for a hyperbolic twist knot at its holonomy representation is expressed by using the cusp shape of the hyperbolic structure of the knot complement. The last part of the paper (Subsection~\ref{ExTb}) deals with some remarks on the behavior of the sequence of non--abelian Reidemeister torsions for twist knots at the holonomy indexed by the number of crossings. The appendix contains concrete examples and tables of the values of the non--abelian torsion for some explicit twist knots.

%%%%%%%%%%%%%%%%%%%%%%%%%%%%%%%%%%%
\section{Twist knots}
\label{section:twist_knots}
%%%%%%%%%%%%%%%%%%%%%%%%%%%%%%%%%%%
\label{TwistKnots}

\begin{notation*}
According to the notation of \cite{HS}, 
the twist knots  are written $J(\pm 2, n)$, 
where $n$ is an integer. The $n$ crossings are right--handed when $n>0$ and left--handed when $n<0$.
%%% figure of J(2, n)
\begin{figure}[htb]\label{fig:J(2,n)}
\begin{center}
\includegraphics{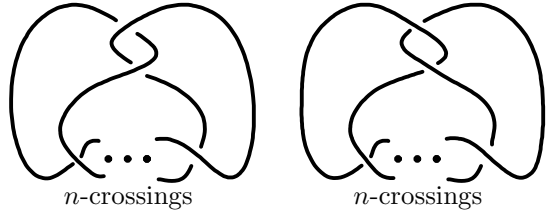}
\end{center}
\caption{The diagrams of $J(2, n)$ and $J(-2,n )$, $n>0$.}
\end{figure}
\end{notation*}

Here are some important facts about twist knots.
\begin{enumerate}
\item By definition, if $n \in \{0, 1\}$ then $J(2,n)$ and $J(-2,-n)$ are the
unknot. Also {observe that} $J(2,-1)$ and $J(-2,1)$ are isotopic to $J(-2,-2)$ and $J(2,2)$
respectively (for general cases, see below). 
{For this reason, in all this paper,} we focus on the knots $J(\pm 2, n)$ with $|n| \geqslant 2$.

\item If we rotate the diagram of $J(2,n)$ by a 90 degrees angle clockwise then we get a diagram of a rational knot in the sense of Conway. 
In rational knot  notation $J(2,n)$, $n>0$, is represented by the continued fraction 
$[n,-2]=\frac{1}{-2 + 1/n}=\frac{-n}{2n-1}.$ 
Therefore in two--bridge knot notation, for $n>0$ we have 
$J(2,n)=b(2n-1,-n)=b(2n-1,n-1).$ 

Similarly, the knot $J(2,-n)$, $n>0$, is represented by the continued fraction $[-n,-2]$, therefore $J(2, -n) = b(2n+1,n+1).$ 
On the other hand, the knot $J(-2,n)$, $n>0$,  is $b(2n+1,n)$ and $J(-2,-n)$, $n>0$,  is $b(2n-1,n)$. 
  \item {Two other important observations are} the following:
  \begin{itemize}
  \item the twist knot $J(2,n)$ is isotopic to $J(-2,n-1)$ (see~\cite{HS}),
  \item and {$J(-2,-n)$ is the mirror image of $J(2,n)$.}
\end{itemize}
  \end{enumerate}
 
As a consequence,  
\emph{we will only consider the twist knots $J(2, n)$, 
where $n$ is an integer such that $|n| \geqslant 2$}. 
From now on, we adopt in the sequel the following terminology: 
a twist knot $J(2, n)$ is said to be \emph{even} (resp. \emph{odd}) 
if $n$ is even (resp. odd).

\begin{example*}
Note that in Rolfsen's table~\cite{Rolfsen}, 
the trefoil knot $3_1 = J(2, 2) = b(3, 1)$, 
the figure eight knot $4_1 = J(2, -2) = b(5,3)$, $5_2 = J(2, 4)$ and $6_1 = J(2, -4)$ etc.
\end{example*}

\begin{notation*}
For a knot $K$ in $S^3$, we let $E_K$ (resp. $\Pi(K)$) 
denote the exterior (resp. the group) of $K$, 
i.e. $E_K = S^3 \setminus N(K)$, 
where $N(K)$ is an open tubular neighborhood of $K$ 
(resp. $\Pi(K) = \pi_1(E_K)$). 
\end{notation*}

\begin{convention*}
Suppose that $S^3$ is oriented. The exterior of a knot is thus oriented and we know that it is bounded 
by a $2$-dimensional torus $T^2$. 
This boundary inherits an orientation by the convention 
\emph{\lq\lq the inward pointing normal vector in the last position\rq\rq}. 
Let $\mathrm{int}(\cdot, \cdot)$ be 
the intersection form on the boundary torus induced by its orientation. 
The peripheral subgroup $\pi_1(T^2)$ is generated by 
the meridian--longitude system $(\mu, \lambda)$ of the knot. 
If we suppose that the knot is oriented, 
then $\mu$ is oriented by the convention that 
the linking number of the knot with $\mu$ is $+1$. 
Next, $\lambda$ is the oriented preferred longitude 
using the rule $\mathrm{int}(\mu, \lambda) = +1$. These orientation conventions will be used in the definition of the {(sign--refined)} twisted Reidemeister torsion.
\end{convention*}
 
 Twist knots live in the more general family of two-bridge knots. 
The group of such a knot admits 
a particularly nice Wirtinger presentation with only two generators and a single relation. 
Such Wirtinger presentations of groups of twist knots  
are given in the two following facts (see for example~\cite{Rolfsen} or~\cite{HS} for a proof). 
We distinguish even and odd cases and suppose that $m \in \ZZ$.
\begin{fact}\label{Fact1}   
    The knot group of $J(2, 2m)$ admits the following presentation:
\begin{equation}\label{Pi1even}
\Pi(J(2, 2m))
=
\langle
x, y \,|\, w^{m} x = y w^{m}
\rangle
\end{equation}
where $w$ is the word $[y, x^{-1}] = yx^{-1}y^{-1}x$.
\end{fact}

\begin{fact}\label{Fact2}
The knot group of $J(2, 2m+1)$ admits the following presentation:
\begin{equation}\label{Pi1odd}
\Pi(J(2, 2m+1))
=
\langle
x, y \,|\, w^{m} x = y w^{m}
\rangle
\end{equation}
where $w$ is the word $[x,y^{-1}] = xy^{-1}x^{-1}y$.
\end{fact}

One can easily describe the peripheral--system $(\mu, \lambda)$ of a twist knot.
It is expressed in the knot group as:
\[
\mu = x \text{ and } \lambda = (\stackrel{\gets}{w})^m w^m,
\]
where we let $\stackrel{\gets}{w}$ denote the word obtained from $w$ by reversing the order of the letters.

\begin{remark}\label{wordtwistknot}
The knot group of a two--bridge knot $K$ admits a {\emph{distinguished}} Wirtinger presentation of the following form:
\begin{equation*}%\label{P:G2B} 
\Pi(K)=\left\langle {x,y \, | \, \Omega x=y\Omega} \right\rangle  \text{ where } \Omega=x^{\epsilon_1}y^{\epsilon_n}x^{\epsilon_2}y^{\epsilon_{n-1}}\dotsb x^{\epsilon_n}y^{\epsilon_1}, \text{ with } \epsilon_i = \pm 1.
\end{equation*}
With the above notation, for $K = J(2, 2m)$, $m \in \ZZ \setminus \{0\}$, the word $\Omega$ is: 
\begin{equation}\label{word}
\Omega_m = \begin{cases}
     w^m & \text{ if } m < 0, \\
      x^{-1} {(\overline{w})}^{m-1} y^{-1} & \text{ if } m > 0.
\end{cases}
\end{equation}
Here $\overline{w} = {(\stackrel{\gets}{w})}^{-1}$, i.e. the word $\overline{w}$ is obtain from $w$ by changing each of its letters by its reverse. Of course this choice is strictly equivalent to presentation~(\ref{Pi1even}). But in a sense, when $m>0$ the word $w^m$ does not give a \lq\lq reduced" relation (some cancelations are possible in $w^mxw^{-m}y^{-1}$) which is not the case for {the word} $\Omega_m$.
\end{remark}

Some more elementary properties of twist knots are discussed in the following remark.
\begin{remark}\label{remarktwist}
\begin{enumerate} 
\item \emph{The knot groups $\Pi(J(2, 2m+1))$ and $\Pi(J(2, -2m))$ are isomorphic} by interchanging $x$ and $y$. {Algebraically, it is thus enough to consider the case of even twist knots (see in particular Remark~\ref{rmsgn}).}
  \item \emph{The genus of a twist knot is $1$} (\cite[p. 99]{Adams}). Thus, the only torus knot which is a twist knot is the trefoil knot $3_1$.
  \item \emph{Twist knots are hyperbolic knots} except  the trefoil knot (see for exam\-ple~\cite{Menasco}).
  
  \item It is well known (see for example~\cite{Rolfsen}) that 
  \emph{the Alexander polynomial of the twist knot $J(2, 2m)$ is given by
  $\Delta_{J(2, 2m)}(t) = mt^2 + (1-2m)t + m.$}
  Moreover, using the mirror image invariance of the Alexander polynomial, one has $\Delta_{J(2, 2m+1)}(t) = \Delta_{J(2, -2m)}(t)$. 
Thus the Alexander polynomial becomes monic if and only if 
$m$ is $\pm 1$. As a consequence,  \emph{the knot $J(2, 2m)$ is not fibered} (since its Alexander polynomial is not monic) except for $m = \pm 1$, that is to say except for the trefoil knot and the figure eight knot, which are known to be fibered knots.
\begin{figure}[!htb]
\begin{center}
\scalebox{.9}{\includegraphics{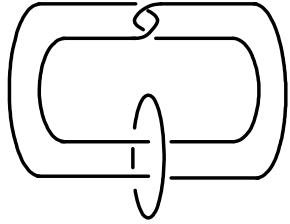}}
\caption{The Whitehead link.}
\label{fig:Whitehead}
\end{center}
\end{figure}
    \item Twist knot exteriors can be obtained by surgery on the trivial component of the Whitehead link $\mathcal{W}$ (see Fig.~\ref{fig:Whitehead}). More precisely, \emph{$E_{J(2, -2m)} = \mathcal{W}({1/m})$ is obtained by a surgery of slope $1/m$} on the trivial component of the Whitehead link $\mathcal{W}$, see~\cite[p.~263]{Rolfsen} for a proof.
As a consequence, twist knots are \emph{all virtually fibered} ({a manifold $M$ is said to be virtually fibered if there is a finite cover of $M$ that is fibered}), see \cite{Leininger}.
\end{enumerate}
\end{remark}

%%%%%%%%%%%%%%%%%%%%%%%%%%%%%%%%%%%%%%%%%%%%%%%%%%%
\section{On the $\SL$-character variety and non--abelian representations}
\label{section:review_character_variety}
%%%%%%%%%%%%%%%%%%%%%%%%%%%%%%%%%%%%%%%%%%%%%%%%%%%
\subsection{Review on the $\SL$-character variety of knot groups}
\label{ModuliK}

Given a finitely genera\-ted group $\pi$, we let $$R(\pi; \SL) = \mathrm{Hom}(\pi; \SL)$$ denote the space of $\SL$-representations of $\pi$. As usual, this space is endowed with the {compact--open topology}. Here $\pi$ is assumed to have the discrete topology and the Lie group $\SL$ is endowed with the usual one. 

A representation $\rho \colon \pi \to \SL$ is called \emph{abelian} if $\rho(\pi)$ is an abelian subgroup of $\SL$.  A representation $\rho$ is called \emph{reducible} if there exists a \emph{proper} subspace $U \subset \C^2$ such that $\rho(g)(U) \subset U$, for all $g \in \pi$. Of course, any abelian representation is reducible (while the converse is false in general). A  
representation is called {\it irreducible} if it is not reducible.

  The group $\SL$ acts on the representation space $R(\pi; \SL)$ by conjugation, but the naive quotient $R(\pi; \SL)/\SL$ is not Hausdorff in general. Following~\cite{CS:1983}, we will focus on the \emph{character variety} $X(\pi) = X(\pi; \SL)$ which is the set of \emph{characters} of $\pi$. Associated to the representation $\rho \in R(\pi; \SL)$, its character $\chi_\rho \colon \pi \to \C$ is defined by $\chi_\rho(g) = \mathrm{tr}(\rho(g))$, where $\mathrm{tr}$ denotes the trace of matrices. 
    In some sense $X(\pi)$ is the \lq\lq algebraic quotient\rq\rq\, of $R(\pi; \SL)$ by the action by conjugation of $\PSL$. It is well known that $R(\pi, \SL)$ and $X(\pi)$ have the structure of {complex algebraic affine varieties} (see~\cite{CS:1983}).

  Let ${R}^\mathrm{irr}(\pi; \SL)$ denote the subset of irreducible representations of $\pi$, and $X^\mathrm{irr}(\pi)$ denote its image under the map ${R}(\pi; \SL) \to X(\pi)$. Note that two irreducible representations of $\pi$ in $\SL$ with the same character are conjugate by an element of $\SL$, see~\cite[Proposition 1.5.2]{CS:1983}. Similarly, we write $X^\mathrm{nab}(\Pi(K))$ for the image of the set $R^{\mathrm{nab}}(\Pi(K))$ of non--abelian representations in $X(\Pi(K))$. 
Note that $X^\mathrm{irr}(\Pi(K)) \subset X^\mathrm{nab}(\Pi(K))$ and observe that this inclusion is strict in general.

%For a knot $K$ in $S^3$, let $\Pi(K)$ denote the knot group. 
%Let $\Pi(K)'$ be the subgroup generated by the commutators of $\Pi(K)$. It is well known that $\Pi(K)/\Pi(K)' \cong H_1(S^3 \setminus K; \ZZ) \cong \ZZ$ is generated by the meridian $\mu$ of $K$. As a consequence, each {abelian} representation of $\Pi({K})$ in $\SL$ is conjugate either to $\varphi_z \colon \Pi({K}) \to \SL$ defined by $\varphi_z(\mu) = \left(\begin{array}{cc}e^z & 0 \\0 & e^{-z}\end{array}\right)$, with $z \in \C$, or to a representation $\rho$ with $\rho(\mu) = \pm \left(\begin{array}{cc}1 & 1 \\0 & 1\end{array}\right)$.  
%
%Let $R^\mathrm{ab}(\Pi(K))$ denote the space of abelian representations, its image under the map $R(\Pi(K); \SL) \to X(\Pi(K))$ is denoted $X^\mathrm{ab}(\Pi(K))$. 

%If $M$ is a $3$-dimensional manifold, then we write 
%$$
%X^\mathrm{ab}(M) = X^\mathrm{ab}(\pi_1(M)), \; 
%X^\mathrm{nab}(M) = X^\mathrm{nab}(\pi_1(M)), \;X^\mathrm{irr}(M) = X^\mathrm{irr}(\pi_1(M)) %\text{ etc.}
%$$

\subsection{Review on  the character varieties of  two--bridge knots}
\label{ReviewCharK}
Here we briefly review Riley's method~\cite{Riley} for describing the non--abelian part of the representation space of two--bridge knot groups.

The knot group of a two--bridge knot $K$ admits a presentation of the following form:
\begin{equation}\label{P:G2B} 
\Pi(K)=\left\langle {x,y \, | \, \Omega x=y \Omega} \right\rangle  \text{ where } \Omega=x^{\epsilon_1}y^{\epsilon_n}x^{\epsilon_2}y^{\epsilon_{n-1}}\dotsb x^{\epsilon_n}y^{\epsilon_1}, \; \epsilon_i = \pm 1.
\end{equation}
  
  We consider the following matrices: 
%$C=  \begin{pmatrix}    t&1\\    0&1  \end{pmatrix},\ D=  \begin{pmatrix}  \hphantom{-t}t&0\\  -tu&1  \end{pmatrix};$
$C_1=
  \begin{pmatrix}
    \sqrt{s}&1\\
    0&{1}/{\sqrt{s}}
  \end{pmatrix},\ 
D_1=
  \begin{pmatrix}
   \hphantom{-}\sqrt{s}&0\\
   -u&{1}/{\sqrt{s}}
  \end{pmatrix};
$
and 
$
C_2=
  \begin{pmatrix}
    \sqrt{s}&{1}/{\sqrt{s}}\\
    0&{1}/{\sqrt{s}}
  \end{pmatrix},\ 
D_2=
  \begin{pmatrix}
    \hphantom{-t}\sqrt{s}&0\\
    -\sqrt{s}u&{1}/{\sqrt{s}}
  \end{pmatrix};
$
where $s \in \CC^*$, $\sqrt{s}$ is some square root, and $u \in \CC$.
\begin{remark}\label{remCD}
Note that $C_2$ and $D_2$ are respectively conjugate to $C_1$ and $D_1$ via $U=\bigl(
   \begin{smallmatrix}
     s^{-1/4}&0\\
     0&s^{1/4}
   \end{smallmatrix}
\bigr)$.
\end{remark}

%\begin{remark}\label{remCD}
%Note that $C$ and $D$ can be obtained by conjugating $tC_1$ and $tD_1$ 
%by the diagonal matrix 
%%$U=\bigl(
 %  \begin{smallmatrix}
 %    t^{-1/2}&0\\
 %    0&t^{1/2}
 %  \end{smallmatrix}
%\bigr)$, and replacing $t^2$ by $t$.
%We also note that $C_2$ and $D_2$ are conjugate to $C_1$ and $D_1$ via $U$.
%\end{remark}
 
%\begin{fact}\label{Cl1}
%If $M_1$, $M_2$ are non--commuting elements in $\SL$ with same traces, then there exists a pair $(t, u)\in \CC^2$ such that $M_1$ and $M_2$ are conjugated to $C_2$ and $D_2$ respectively.
%\end{fact}

 %Combining Fact~\ref{Cl1} and~Remark~\ref{remCD}, we obtain:\begin{claim}If $M_1$, $M_2$ are non-commuting elements in $\SL$ with same traces, then $M_1$ and $M_2$ are simultaneously conjugated to $C_1$ and $D_1$ respectively.\end{claim}

For a two--bridge knot $K$, 
$x$ and $y$ are conjugate elements in $\Pi(K)$ and represent meridians of the knot,
 therefore $\rho(x)$ and $\rho(y)$ have same traces. 
If $\rho(x)$ and $\rho(y)$ do not commute, i.e. if $\rho$ is non--abelian, 
then up to a conjugation 
one can assume that $\rho(x)$ and $\rho(y)$ are the matrices $C_1$ and $D_1$ 
respectively. With this notation, we have the following result due to Riley~\cite{Riley}.

%\begin{proposition}The homomorphism $\rho \colon \Pi(K) \to \mathrm{GL}_2(\CC)$ defined by  $\rho(x)=C$ and $\rho(y)=D$ is a non--abelian representation of $\Pi(K)$  if and only if the pair $(t,u) \in \CC^2$ satisfies the following equation \begin{equation}\label{eqn:def_Riley_poly}  w_{1, 1}+(1-t)w_{1, 2}=0\end{equation} where $W = \rho(\Omega) =\left(   \begin{array}{cc}     w_{1, 1} & w_{1, 2} \\     w_{2, 1} & w_{2, 2}   \end{array}  \right)$. 

%Conversely, every non--abelian representation is conjugated to a representation satisfying Equation~(\ref{eqn:def_Riley_poly}).\end{proposition}

%\begin{proof} A direct matrix computation shows that the requirement $W C=D W$ is equivalent to the following two equations: 
%$$ w_{1, 1}=(t-1)w_{1, 2} \text{ and } w_{2, 1}=-tuw_{1, 2}.$$ 
%The second equation is just a consequence of the fact that $W$ is palindromic and therefore is not really a requirement at all. Indeed, using $W^T$ to denote the transpose of the matrix $W$, we have $W^T={D^T}^{\epsilon_1}{C^T}^{\epsilon_n}{D^T}^{\epsilon_2}{C^T}^{\epsilon_{n-1}}\dotsb {D^T}^{\epsilon_n}{C^T}^{\epsilon_1}$. Now $C^T=VDV^{-1}$ and $D^T=VCV^{-1}$ where $V=\bigl(\begin{smallmatrix}(-tu)^{1/2}&0\\0&(-tu)^{-1/2}\end{smallmatrix}\bigr)$. This provides $V^{-1}W^TV=W$, which immediately gives the second equation above.\end{proof}

%\begin{notation*}We let $\phi_K(t, u) = w_{1, 1}+(1-t)w_{1, 2}$ denote the left hand side of Equation $(\ref{eqn:def_Riley_poly})$ and call it the \emph{Riley polynomial} of $K$. \end{notation*} %The same proof provides a sometimes more convenient result:
\begin{proposition} 
The homomorphism $\rho \colon \Pi(K) \to \mathrm{SL}_2(\CC)$ 
defined by $\rho(x)=C_1$ and $\rho(y)=D_1$ 
is a non--abelian representation of $\Pi(K)$ 
if and only if the pair $(s,u) \in \CC^2$ satisfies the following equation: 
\begin{equation}\label{eq:R}
  w_{1, 1}+({1}/{\sqrt{s}}-\sqrt{s})w_{1, 2}=0
\end{equation}
where $W=\rho(\Omega)$. 

Conversely, every non--abelian representation is conjugated to a representation satisfying Eq.~(\ref{eq:R}).
\end{proposition}

\begin{remark}
Similarly, if $\rho(x)=C_2$ and $\rho(y)=D_2$, then Riley's equation is 
\begin{equation}\label{eq:R_alternative}
w_{1, 1}+(1-s)w_{1, 2}=0.
\end{equation}
\end{remark}

\begin{notation}
We let $\phi_K(s, u) = w_{1, 1}+(1-s)w_{1, 2}$ denote the left hand side of
Eq. $(\ref{eq:R_alternative})$ 
and call it the \emph{Riley polynomial} of $K$.
\end{notation}

\subsection{The holonomy representation of a hyperbolic twist knot}
\label{sechyp}
\subsubsection{Some generalities}
    It is well known that the complete hyperbolic structure of a hyperbolic knot complement determines a unique discrete faithful representation of the knot group in $\mathrm{PSL}_2(\CC)$, called the \emph{holonomy representation}. It is proved~\cite[Proposition 1.6.1]{HandBookGT} that such a representation lifts to $\SL$ and determines two irreducible representations in $\SL$.
    
    The trace of the peripheral--system at the holonomy is $\pm 2$ because their images by the holonomy are parabolic matrices. More precisely, Calegari~\cite{Calegari} proved that the trace of the longitude at the holonomy is always $-2$ and the trace of the meridian at the holonomy is $\pm 2$, depending on the choice of the lift. 
We summarize all this in the following important fact.

\begin{fact}\label{fact:Calegari}
Let $\rho_0$ be one of the two lifts of the discrete and faithful representation associated to the complete hyperbolic structure of a hyperbolic knot $K$ and let {$T^2 = \partial E_K$} denote the boundary of the knot exterior. The restriction of $\rho_0$ to $\pi_1(T^2)$ is conjugate to the parabolic representation such that
\[
\mu \mapsto 
\pm \left(
    \begin{array}{cc}
      1 & 1 \\
      0 & 1
    \end{array}
    \right), 
\quad 
\lambda \mapsto 
        \left(
        \begin{array}{cc}
          -1 & \hphantom{-}\cusp \\
          \hphantom{-}0 & -1
        \end{array}
        \right).
\]
Here $\cusp = \cusp(\lambda, \mu) \in \CC$ is called the cusp shape of $K$.
\end{fact}

\begin{remark}
The universal cover of the exterior of a hyperbolic knot is the hyperbolic $3$-space $\mathbb{H}^3$. The cusp shape can be seen as the ratio of the translations of the parabolic isometries of $\mathbb{H}^3$ induced by projections to $\mathrm{PSL}_2(\CC)$. Of course, the cusp shape $\cusp = \cusp(\lambda, \mu)$ depends on the choice of the basis $(\mu, \lambda)$ for $\pi_1 (T^2)$. A change in the basis of $\pi_1 (T^2)$ shifts $\cusp$ by an integral M\"obius transformation. 
\end{remark}

\subsubsection{Holonomy representations of twist knots}

This subsection is concerned with the $\SL$-representations which are lifts of the holonomy representation in the special case of (hyperbolic) twist knots. In particular, we want to specify the images, up to conjugation, of the group generators $x$ and $y$ (see the group presentation~(\ref{P:G2B})).

\begin{lemma}\label{L8}
Let $K$ be a hyperbolic two--bridge knot and 
suppose that its knot group admits the following presentation 
$\Pi(K)=\langle x,y \, | \, \Omega x=y \Omega\rangle$. 
If $\rho_0$ denotes a lift in $\SL$ of the holonomy representation, 
then $\rho_0$ is given, up to conjugation, {by}
\[
x \mapsto 
 \pm \left(
\begin{array}{cc}
 1 & 1 \\
 0 & 1 
\end{array}
\right)
, \quad
y \mapsto 
 \pm\left(
\begin{array}{cc}
 \hphantom{-}1 & 0 \\
 - u & 1 
\end{array}
\right),
\]
where $u$ is a root of Riley's equation $\phi_K(1, u) = 0$ of $K$.
\end{lemma}

\begin{proof}
It follows from Fact~\ref{fact:Calegari} that
each lift of the holonomy representation maps the meridian to 
$\pm \left(
\begin{smallmatrix}
1 & 1\\
0 & 1
\end{smallmatrix}
\right).$
It is known that the lifts of the holonomy representation are 
irreducible $\SL$-representations, 
in particular,  non--abelian ones.
Hence we can construct the $\SL$-representations which are conjugate to 
the lifts of the holonomy representation 
by using roots of Riley's equation.
%Since both generators $x$ and $y$ are meridian, the corresponding root $(s, u)$ of Riley's polynomial satisfies that $\sqrt{s}+1/\sqrt{s}=\pm 2$. Therefore a lift of the holonomy representation corresponds to a root of $(1, u)$ of Riley's polynomial.
Using Section \ref{ReviewCharK}, if $x$ is sent to $\pm \left(
\begin{smallmatrix}
1 & 1\\
0 & 1
\end{smallmatrix}
\right)$ then $y$ is sent to $\pm\left(
\begin{smallmatrix}
 \hphantom{-}1 & 0 \\
 - u & 1 
\end{smallmatrix}\right)$, where $u$ is a root of Riley's equation $\phi_K(1,u) = 0$.
\end{proof}

\begin{notation*}
If we let $A$ be an element of $\SL$, 
then the adjoint actions of $A$ and $-A$ are same.
So, we use the $\SL$-representation such that 
$$
x \mapsto 
\left(
\begin{array}{cc}
  1 & 1 \\
  0 & 1 
\end{array}
\right)
, \quad
y \mapsto 
\left(
\begin{array}{cc}
 \hphantom{-}1 & 0 \\
-u &  1 
\end{array}
\right) \quad (\text{where } u \text{ is such that } \phi_K(1, u) = 0)
$$
as a lift of the holonomy representation and we improperly call it the \emph{holonomy representation}. 
\end{notation*}

\subsection{On parabolic representations of twist knot groups}
\label{ParRep}

In this subsection, we are interested in the \emph{parabolic representations} of (hyperbolic) two--bridge knot groups and especially twist knot groups. The holonomy representation is one of them. Lemma~\ref{L8} characterizes the holonomy representation algebraically and says that it corresponds to a root of Riley's equation $\phi_K(1,u) = 0$. A natural and interesting question is the following: which roots of Riley's equation $\phi_K(1,u) = 0$ correspond to the holonomy representation? Here we will give a geometric characterization of such roots {in the case of two--bridge knots}.

\subsubsection{Remarks}
%\begin{remark}
We begin this section by some elementary but important remarks on the roots of Riley's equation $\phi_K(1,u) = 0$ corresponding to holonomy  representations.
\begin{enumerate}
  \item One can first notice that such roots are necessarily \emph{complex numbers which are not real}, because the discrete and faithful representation is irreducible and \emph{not} conjugate to a real representation (i.e. a representation such that the image of each element is a matrix with real entries).
  \item One can also observe that holonomy representations correspond to \emph{a pair of complex conjugate roots} of  Riley's equation $\phi_K(1,u) = 0$ as it is easy to see.
\end{enumerate}
%\end{remark}

\subsubsection{Generalities: the case of two--bridge knots}
Let $K$ be a hyperbolic two--bridge knot.
Suppose that a presentation of the knot group $\Pi(K)$ is given as in Eq.~(\ref{P:G2B}) by
\[
\Pi(K) = \langle {x, y \, |\, \Omega x = y\Omega} \rangle, \text{ where } \Omega \text{ is a word in } x,y. 
\]
The longitude of $K$ is of the form: $\lambda = \stackrel{\gets}{\Omega}\Omega \, x^n$. Here $n$ is an integer such that 
the sum of the exponents in the word $\lambda$ is $0$ and we repeat that $\stackrel{\gets}{\Omega}$ denotes the word obtained from $\Omega$ by reversing the order of the letters.

Let $\rho \colon \Pi(K) \to \SL$ be a representation such that: 
\[
x \mapsto 
\pm \left(
\begin{array}{cc}
1 &  1 \\
  0          &  1 
\end{array}
\right)
, \quad
y \mapsto 
\pm \left(
\begin{array}{cc}
 1 & 0 \\
  - u &  1 
\end{array}
\right)
\]
%Following the results in~\cite{Riley}, if $\rho$ is an $\SL$-representation of the knot group $\Pi(K)$, then 
where $u$ is necessarily a root of Riley's equation $\phi_K(1,u)=0$. 
Suppose that 
\[
\rho(\Omega) = \left(\begin{array}{cc}w_{1,1} & w_{1,2} \\w_{2,1} & w_{2,2}\end{array}\right)
\]
where $w_{i,j}$ is a polynomial in $u$ for all $i, j \in \{1,2\}$.

Riley's method gives us the following identities (see Section~\ref{ReviewCharK}):
\[
w_{1,1} = 0 \text{ and } u w_{1,2} + w_{2,1} = 0. 
\]
Thus
\[
\rho(\Omega) = \left(\begin{array}{cc}0 & w_{1,2} \\ -u w_{1,2} & w_{2,2}\end{array}\right).
\]
The fact that $\rho(w) \in \SL$ further gives the following equation:
\begin{equation}\label{UW}
u w_{1,2}^2 = 1.
\end{equation}

The crucial point for computing $\rho(\lambda)$ is to express $\rho(\stackrel{\gets}{\Omega})$ with the help of $\rho(\Omega)$.  
Consider the diagonal matrix:
\[
\mathbf{i} = \left(\begin{array}{cc}i & 0 \\0 & -i\end{array}\right) \in \SL,
\]
where $i$ stands for a square root of $-1$. Let $Ad$ denote the adjoint representation of the Lie group $\SL$. Then the following identities hold:
\[
\rho(x^{-1}) = Ad_{\mathbf{i}} (\rho(x)), \; 
\rho(y^{-1}) = Ad_{\mathbf{i}} (\rho(y))  
\text{ and } 
\rho(\Omega^{-1}) =  Ad_{\mathbf{i}} (\rho(\stackrel{\gets}{\Omega})).
\]
Thus, we have
\[
 \rho(\stackrel{\gets}{\Omega}) 
= Ad_{\mathbf{i}} (\rho({\Omega}^{-1})) 
= \left(
  \begin{array}{cc}
    w_{2,2} & w_{1,2} \\
    w_{2,1} & w_{1,1}
  \end{array}
  \right).
\]
Next, a direct computation gives:
\begin{equation}\label{long}
\rho(\lambda) = \rho(\stackrel{\gets}{\Omega}) \rho(\Omega) \rho(x)^n =\left(\begin{array}{cc}-u w^2_{1,2} & -n u w^2_{1,2} + 2w_{1,2}w_{2,2} \\0 & -u w^2_{1,2}\end{array}\right).
\end{equation}
Combining Eqs.~(\ref{UW}) and~(\ref{long}), we obtain
\begin{equation}\label{long2}
\rho(\lambda) = \left(\begin{array}{cc}-1 & -n + 2w_{1,2}w_{2,2} \\0 & -1\end{array}\right).
\end{equation}
And we conclude that the \emph{cusp shape of $K$} is 
\begin{equation}\label{cuspshape_bridge}
\cusp = n -  2w_{1,2}w_{2,2}.
\end{equation}

\begin{remark}
In particular, Eq.~(\ref{long2}) gives us, by an elementary and direct computation, Calegari's result~\cite{Calegari}: $\mathrm{tr}\, \rho_0(\lambda) = - 2$ for the discrete faithful representation $\rho_0$ associated to the complete hyperbolic structure of the exterior of a (hyperbolic) two--bridge knot.
\end{remark}

\subsubsection{The special case of twist knots}
In the case of hyperbolic twist knots, we can further estimate $w_{i,j}$ in Eq.~(\ref{cuspshape_bridge}). In fact, we only consider the case where $K = J(2, 2m)$ in what follows. 
The group $\Pi(K)$ of such a knot has the following presentation:
\[
\Pi(K) = \langle {x, y \, |\, w^m x = y w^m} \rangle,
\]
where $w$ is the commutator $[y,x^{-1}]$ (see Fact~\ref{Fact1}, Section~\ref{TwistKnots}). A direct computation of the commutator $\rho(w) = [\rho(y),\rho(x)^{-1}]$ gives:
\[
W= \rho(w) = \left(
\begin{array}{cc}
  1-u  & -u \\
  u^2  & u^2+u+1
\end{array}
\right) \quad (\text{where } u \text{ is such that } \phi_K(1,u) = 0).
\]

Using the Cayley--Hamilton identity, it is easy to obtain the following recursive formula for the powers of the matrix $W$:
\begin{equation}\label{CH}
W^k - (u^2+2) W^{k-1} + W^{k-2} = 0, \; k \geqslant 2.
\end{equation}
Eq.~(\ref{CH}) implies, by induction on $k$,
\[
w_{2,2} = - (u+2) w_{1,2}.
\]
Since $n = 0$ and $u w_{1,2}^2 = 1$, the \emph{cusp shape of the twist knot $K$} is:
\[
\cusp = n + (2u+4)w_{1,2}^2 = \frac{2u+4}{u}.
\]

    In other words, the root $u_0$ of Riley's equation $\phi_K(1,u) = 0$ corresponding to the holonomy representation satisfies the following equation:
\begin{equation}\label{cuspshape}
u_0 = \frac{4}{\cusp - 2},
\end{equation}
where $\cusp$ is the cusp shape of the knot exterior.

\begin{remark}
Eq.~(\ref{cuspshape}) gives a geometric characterization of the (pair of complex conjugate) roots of Riley's equation $\phi_K(1,u) = 0$ associated to the holonomy representation in terms of the cusp shape, a geometric quantity associated to each cusped hyperbolic $3$-dimensional manifold.
\end{remark}

\subsection{The character varieties of twist knots: a recursive description}

T. Le~\cite{Le} gives a recursive description of the $\SL$-character variety of two--bridge knots and apply it to obtain an \emph{explicit} description of the $\SL$-character variety of torus knots. 
Here we apply his method to obtain an explicit recurrent description of the $\SL$-character variety of twist knots. This is the main difference with the description given by Hoste and Shanahan~\cite{HS}.

Let $n = 2m$ or $2m+1$, recall that 
$\Pi(J(2,n))=\langle x,y \,|\, \Omega_m x=y \Omega_m \rangle,$
where $\Omega_m$ is a word in $x, y$ (see Facts~\ref{Fact1}--\ref{Fact2}).

\begin{notation*}
Let $\gamma \in \Pi(J(2,n))$. Following a notation introduced in~\cite{DehnSurgery}, we let 
$$I_\gamma : X(\Pi(J(2,n))) \to \CC$$ 
be the trace--function defined by $I_\gamma : \rho \mapsto \tr(\rho(\gamma))$.
\end{notation*}

Let $a=I_x$, $b=I_{xy}$ and recall the following useful formulas for $A,B,C \in \SL$:
\begin{equation}\label{Tr1}
\mathrm{tr}(A^{-1}) = \mathrm{tr}(A) \text{ and } \mathrm{tr}(AB) = \mathrm{tr}(BA),
\end{equation}
\begin{equation}\label{Tr2}
\mathrm{tr}(AB) + \mathrm{tr}(A^{-1}B) = \mathrm{tr}(A)\mathrm{tr}(B),
\end{equation}
\begin{equation}\label{Tr3}
\mathrm{tr}(ABA^{-1}B^{-1}) = -2 - \mathrm{tr}(A)\mathrm{tr}(B)\mathrm{tr}(AB) + (\mathrm{tr}(A))^2 +(\mathrm{tr}(B))^2+(\mathrm{tr}(AB))^2.
\end{equation}

As $x$ and $y$ are conjugate elements in $\Pi(J(2,n))$, we have $I_y = a = I_x$.
If $\gamma$ is a word in the letters $x$ and $y$, then $I_\gamma$ can always be expressed as a polynomial function  in $a$ and $b$. For example, combining the usual Formulas~(\ref{Tr1}),~(\ref{Tr2}) and~(\ref{Tr3}), one can easily observe that for $w=[y, x^{-1}]=yx^{-1}y^{-1}x$, we have
\begin{equation}\label{I_w}
I_w = -2 - a^2b + 2a^2+b^2.
\end{equation}

The character variety of $\Pi(J(2,n))$ is thus parametrized by $a$ and $b$. Here is a practical description of it:
\begin{enumerate}
  \item We first consider the abelian part of the character variety. It is easy to see that the equation $a^2 - b - 2 = 0$ determines the abelian part of the character variety of any knot group.
  \item Next, consider the non-abelian part of the character variety of $\Pi(J(2,n))$, suppose that the length of the word $\Omega_m$ is $2k+2$ (we know that the length of $\Omega_m$ is even). According to \cite[Theorem 3.3.1]{Le}, the non--abelian part of the character variety of $\Pi(J(2,n))$, {for} $n=2m$ or $2m+1$, is determined by the polynomial equation:
\begin{equation*}\label{EQ:CharVar}
{\mathbf\Phi}_{m}(a,b)=0,
\end{equation*}
where
$${\mathbf\Phi}_{m}(a,b)=I_{\Omega_m}-I_{\Omega_m'}+\dotsb+(-1)^{k}I_{\Omega_m^{(k)}}+(-1)^{k+1}.$$ Here we adopt the following notation: if  $\Lambda$ is a word then $\Lambda'$ denotes the word obtained from $\Lambda$ by deleting the two end letters.  
\end{enumerate}

Let us give the two simplest examples to illustrate this general result and find again some well--known facts.

\begin{example}\label{Ex1}
The trefoil knot $3_1$ is the twist knot $J(2,2)$. With the above notation, one has $\Omega_1 = x^{-1}y^{-1}$. Thus applying { Le's method}, the non--abelian part of the character variety is given by the polynomial equation:
\[
{\mathbf\Phi}_{1}(a,b)=I_{x^{-1}y^{-1}} - 1 =0,
\]
which reduces to
\begin{equation*}\label{NabTrefoil}
b = 1.
\end{equation*}
\end{example}

\begin{example}\label{Ex2}
The figure eight knot $4_1$ is the twist knot $J(2,-2)$. With the above notation, one has $\Omega_{-1} = x^{-1}yxy^{-1}$. Thus, the non--abelian part of the character variety of the group of the figure eight knot is given by the polynomial equation:
\[
{\mathbf\Phi}_{-1}(a,b)=I_{x^{-1}yxy^{-1}} - I_{yx} + 1 = 0,
\]
which reduces, using Equation~(\ref{I_w}), to:
\begin{equation*}\label{Nab8}
2a^2+b^2-a^2b - b = 1.
\end{equation*}
\end{example}

Now, we turn back to the general case and only consider the twist knot $J(2, 2m)$ (see Item (2) of Remark~\ref{remarktwist}). Recall that (see Remark~\ref{wordtwistknot}):
\begin{equation}\label{eqWOmega}
\Pi(J(2, 2m)) = \langle x, y \; |\; \Omega_m x = y \Omega_m \rangle, \text{ where } \Omega_m = \begin{cases}
     w^m & \text{ if } m < 0, \\
      x^{-1} {(\overline{w})}^{m-1} y^{-1} & \text{ if } m > 0.
      \end{cases}
\end{equation}
Here $w=[y,x^{-1}]$ and $\overline{w} = [y^{-1}, x]$ and observe that the length of the word $\Omega_m$ is $4m$, if $m < 0$, and $4m-2$, if $m>0$. 

Our method is based on the fact that the word $\Omega_m$ in the distinguished Wirtinger presentation~(\ref{eqWOmega}) of $\Pi(J(2, 2m))$ presents a particularly nice \lq\lq periodic" property. This property is discussed in the following obvious claim.

\begin{claim}\label{claimOmega}
For $m \in \ZZ^*$, we have
\[
    \Omega_{m}^{(4)} = \begin{cases}
    w^{m+2}  & \text{ if } m\leqslant -2,\\
     x^{-1} {(\overline{w})}^{m-3}y^{-1} & \text{ if } m \geqslant 3.
\end{cases}
\]
\end{claim}

Based on Claim~\ref{claimOmega}, for $m \geqslant 0$, we adopt the following notation:
$$S^+_m=I_{\Omega_{m+2}}, \; T^+_m=I_{{(\Omega_{m+2})}'},\; U^+_m=I_{{(\Omega_{m+2})}''}  \text{ and } V^+_m=I_{{(\Omega_{m+2})}'''},$$
and similarly, for $m\leqslant 0$,
$$S^-_{m}=I_{\Omega_{m-2}}, \; T^-_{m}=I_{{(\Omega_{m-2})}'},\; U^-_{m}=I_{{(\Omega_{m-2})}''}  \text{ and } V^-_{m}=I_{{(\Omega_{m-2})}'''}.$$

The Cayley--Hamilton identity applied to the matrix $A^2 \in \SL$ gives$$A^m \left(A^{4}-(\tr A^2)A^{2}+\I\right)=0.$$ 
Write $c(a,b)=I_{\Omega_2}$; thus for $m\geqslant 0$, we have
$$S^+_{m+4} - c(a,b) S^+_{m+2} + S^+_m =0$$ and same relations for $T^+_m$, $U^+_m$ and $V^+_m$ hold. Similarly, we have
$$S^-_{m-4} - c(a,b) S^-_{m-2} + S^-_m =0$$ and same relations for $T^-_m$, $U^-_m$ and $V^-_m$ also hold.

If we write $R^{\pm}_m=S^{\pm}_m-T^{\pm}_m+U^{\pm}_m-V^{\pm}_m$ for $m \in \ZZ$, then above computations can be summarized in the following claim which give us a recursive relation for ${(R^{\pm}_m)}_{m \in \ZZ}$.

\begin{claim} \label{R}
The sequence of polynomials $\left( {R^{\pm}_m(a,b)} \right)_{m \in \ZZ}$  satisfies the following recursive relation:
\begin{equation}\label{R+-}
R^{\pm}_{m\pm 4} - c(a,b) R^{\pm}_{m \pm 2} + R^{\pm}_m =0.
\end{equation}
%\begin{enumerate}
%  \item For $m\geqslant 0$,
%\begin{equation}\label{R+}
%R^+_{m+4} - c(a,b) R^+_{m+2} + R^+_m =0.
%\end{equation}
%  \item For $m \leqslant 0$,
%\begin{equation}\label{R-}
%R^-_{m-4} - c(a,b) R^-_{m-2} + R^-_m =0.
%\end{equation}
%\end{enumerate}
\end{claim}

In Equation~(\ref{R+-}),  using Formula~(\ref{Tr2}) and Equation~(\ref{I_w}), we have:
\[
c(a,b) = I_{w^2} = ({I_w})^2 - 2 = 2+4a^2b-8a^2-4b^2+a^4b^2-4a^4b-2a^2b^3+4a^4+4a^2b^2+b^4.
\]

Let $v$ be such that:
\begin{equation}\label{u/c}
v + v^{-1}=c(a,b).
\end{equation}

For $m \in \ZZ^*$, we distinguish four cases to derive helpful formulas for ${\mathbf\Phi}_m$ in the case of twist knots. 
\begin{itemize}
  \item {\it Case 1:} $m > 0$ is even. 

Let $m=2l$, with $l > 0$, and set $r_i=R_{2i}$. Then 
$$r_{i+2}=c(a,b)r_{i+1}-r_i, \text{ for } i\geqslant 0.$$
 As we have supposed that $v+v^{-1}=c(a,b)$ and following a standard argument in combinatorics (see e.g. \cite[p. 322]{Merris:2003}), we have the general formula (which can also be proved by induction) $r_i=Mv^i+Nv^{-i}$, where $M$ and $N$ are determined by the initial conditions:
\[
\begin{cases}
      r_0=R^+_0=M^+_0+N^+_0  \\
      r_1=R^+_2=M^+_0v+N^+_0v^{-1} 
\end{cases}
\] 
Further observe that:
\[
r_0 = I_{\Omega_2} - I_{\Omega'_2} + I_{\Omega''_2}- I_{\Omega'''_2} = I_{x^{-1}\overline{w}y^{-1}} - I_w + b - 1.
\]
So, we have
\begin{equation*}\begin{split}{\mathbf\Phi}_{m+2}&=R^+_m+R^+_{m-2}+\dotsb+R^+_0\\
&= r_l + \dotsb + r_0\\
&=\sum_{i=0}^l (M^+_0v^i+N^+_0v^{-i})\\&=M^+_0\frac{v^{l+1}-1}{v-1}+N^+_0\frac{v^{-l-1}-1}{v^{-1}-1}.
\end{split}\end{equation*}

  \item {\it Case 2:} $m < 0$ is even. 

Let $m= - 2l$, with $l > 0$, and set $r_i=R^-_{-2i},\ i\geqslant 0$. Similar to the first case, 
$$r_{i+2}=c(a,b)r_{i+1}-r_i, \text{ for } i\geqslant 0.$$ 
The initial conditions are
\[
\begin{cases}
      r_0=R^-_0=M^-_0+N^-_0  \\
      r_1=R^-_{-2}=M^-_0v+N^-_0v^{-1}
\end{cases}
\] 
Further observe that:
\[
r_0 = I_{\Omega_{-2}} - I_{\Omega'_{-2}} + I_{\Omega''_{-2}}- I_{\Omega'''_{-2}} = I_{w^{2}} - I_{y (\overline{w})^{-1}x} + I_w - b.
\]

Thus, we have
\begin{equation*}\begin{split}{\mathbf\Phi}_{m-2}&=R^-_m+R^-_{m+2}+\dotsb+R^-_{0}+1\\
&= r_l + \dotsb + r_0+1\\
&=\sum_{i=0}^l (M^-_0v^i+N^-_0v^{-i})+1\\&=M^-_0\frac{v^{l+1}-1}{v-1}+N^-_0\frac{v^{-l-1}-1}{v^{-1}-1}+1.
\end{split}\end{equation*}

  \item {\it Case 3:} $m > 0$ is odd. 

Let $m=2l+1$, with $l>0$, and set $r_i=R^+_{2i+1},\ i\geqslant 0$. Similar to the first case, 
$$r_{i+2}=c(a,b)r_{i+1}-r_i, \text{ for } i\geqslant 0.$$ 
The initial conditions are
\[
\begin{cases}
      r_0=R^+_1=M^+_1+N^+_1  \\
      r_1=R^+_3=M^+_1v+N^+_1v^{-1}
\end{cases}
\] 
Further observe that:
\[
r_0 = I_{\Omega_3} - I_{\Omega'_3} + I_{\Omega''_3}- I_{\Omega'''_3} = I_{x^{-1}(\overline{w})^{2}y^{-1}} - I_{(\overline{w})^{2}} + I_{xwy} - I_w.
\]
Thus, 
\begin{equation*}\begin{split}{\mathbf\Phi}_{m+2}&=R^+_m+R^+_{m-2}+\dotsb+R^+_1+b-1\\ 
&= r_l + \dotsb + r_0 + b - 1\\
&=\sum_{i=0}^l (M^+_1v^i+N^+_1v^{-i})+b-1\\
&=M^+_1\frac{v^{l+1}-1}{v-1}+N^+_1\frac{v^{-l-1}-1}{v^{-1}-1}+b-1.
\end{split}\end{equation*}

\item {\it Case 4:} $m < 0$ is odd. 

Let $m=-2l-1$, with $l>0$, and set $r_i=R^-_{-2i-1},\ i\geqslant 0$. Similar to the first case, 
$$r_{i+2}=c(a,b)r_{i+1}-r_i, \text{ for } i\geqslant 0.$$ 
The initial conditions are
\[
\begin{cases}
      r_0=R^-_{-1}=M^-_1+N^-_1  \\
      r_1=R^-_{-3}=M^-_1v+N^-_1v^{-1}
\end{cases}
\] 
Further observe that:
\[
r_0 = I_{\Omega_{-3}} - I_{\Omega'_{-3}} + I_{\Omega''_{-3}} - I_{\Omega'''_{-3}} = I_{w^3} - I_{y(\overline{w})^{-2}x} + I_{w^2} - I_{y^{-1}w^{-1}x^{-1}}.
\]

Similarly to the previous case, we have:
\begin{equation*}\begin{split}{\mathbf\Phi}_{m-2}&=R^-_m+R^-_{m+2}+\dotsb +R^-_{-1} + I_{w} - b + 1\\ 
&= r_l + \dotsb + r_0 +I_{w} - b + 1\\
&=\sum_{i=0}^l (M^-_1v^i+N^-_1v^{-i}) +I_{w} - b + 1\\
&=M^-_1\frac{v^{l+1}-1}{v-1}+N^-_1\frac{v^{-l-1}-1}{v^{-1}-1} - a^2b + 2a^2 + b^2 - b - 1.
\end{split}\end{equation*}
\end{itemize}

If we adopt the following notation:
\begin{equation}
\label{MuNu}
M^\pm_j(v, l) = M^\pm_j\frac{v^{l+1}-1}{v-1} \text{ and } N^\pm_j(v, l) = N^\pm_j\frac{v^{-l-1}-1}{v^{-1}-1}, \quad j = 0, 1,
\end{equation}
then we summarize our computations in the following proposition.
\begin{proposition}
The polynomial equation which describes the character variety of the group of the twist knot $J(2,n)$, {for} $n = 2m$ or $2m+1$, is given by  ${\mathbf\Phi}_{m}(a,b) = 0$ where the sequence ${\left( {\mathbf\Phi}_{m}(a,b) \right)}_{m \in \ZZ}$ is recursively defined as follows:
$${\mathbf\Phi}_0(a,b) = 1, \quad {\mathbf\Phi}_{1}(a,b) = b-1, \quad {\mathbf\Phi}_{-1}(a,b) = - a^2b + 2a^2 + b^2 - b - 1,$$ and for $m>1$
\begin{equation}\label{CharVarTK}
{\mathbf\Phi}_{m + 2}(a,b) = \begin{cases}
     M^+_0(v, l) + N^+_0(v, l) & \text{ if } m = 2l \text{ is even}, \\
     M^+_1(v, l) + N^+_1(v, l)+ b - 1 & \text{ if } m=2l+1 \text{ is odd},
\end{cases}
\end{equation}
and
\begin{equation}\label{CharVarTK'}
{\mathbf\Phi}_{-m - 2}(a,b) = \begin{cases}
     M^-_0(v, l) + N^-_0(v, l) + 1 & \text{ if } m = 2l \text{ is even}, \\
     M^-_1(v, l) + N^-_1(v, l) - a^2b + 2a^2 + b^2 - b - 1& \text{ if } m=2l+1 \text{ is odd}.
\end{cases}
\end{equation}
Here $v$ is defined in Equation~(\ref{u/c}) and $M^\pm_k(v, l), N^\pm_k(v, l)$ in Equation~(\ref{MuNu}).
\end{proposition}

\begin{remark}
Observe that in Equations~(\ref{CharVarTK}) and~(\ref{CharVarTK'}), the part $M^\pm_k(v, l) + N^\pm_k(v, l)$ is the \lq\lq recursive" part. In Equation~(\ref{CharVarTK}) the part $b-1$ corresponds to ${\mathbf\Phi}_{1}(a,b)$, and in Equation~(\ref{CharVarTK'}) the part $- a^2b + 2a^2 + b^2 - b - 1$ corresponds to ${\mathbf\Phi}_{-1}(a,b)$, see Examples~\ref{Ex1} and~\ref{Ex2}.
\end{remark}

%%%%%%%%%%%%%%%%%%%%%%%%%%%%%%%%%%%
\section{Review on the non--abelian Reidemeister torsion and twisted polynomial torsion}
\label{section:review_torsion}
%%%%%%%%%%%%%%%%%%%%%%%%%%%%%%%%%%%
\subsection{Preliminaries: the sign-determined torsion of a CW-complex}

We review the basic notions and results about the sign--determined Reidemeister torsion introduced by Turaev which are needed in this paper. Details can be found in Milnor's survey~\cite{Milnor:1966} and in Turaev's monograph~\cite{Turaev:2002}.

\subsubsection*{Torsion of a chain complex}
Let $C_* = (\xymatrix@1@-.5pc{0 \ar[r] & C_n \ar[r]^-{d_n} & C_{n-1} \ar[r]^-{d_{n-1}} & \cdots \ar[r]^-{d_1} & C_0 \ar[r] & 0})$ be a chain complex of finite dimensional vector spaces over $\mathbb{C}$. Choose  a basis $\mathbf{c}^i$ of $C_i$ and  a basis $\mathbf{h}^i$ of the $i$-th homology group $H_i = H_i(C_*)$. The torsion of $C_*$ with respect to these choices of bases is defined as follows.

Let $\mathbf{b}^i$ be a sequence of vectors in $C_{i}$ such that $d_{i}(\mathbf{b}^i)$ is a basis of $B_{i-1}= \im(d_{i} \colon C_{i} \to C_{i-1})$ and let $\widetilde{\mathbf{h}}^i$ denote a lift of $\mathbf{h}^i$ in $Z_i = \ker(d_{i} \colon C_i \to C_{i-1})$. The set of vectors $d_{i+1}(\mathbf{b}^{i+1})\widetilde{\mathbf{h}}^i\mathbf{b}^i$ is a basis of $C_i$. Let $[d_{i+1}(\mathbf{b}^{i+1})\widetilde{\mathbf{h}}^i\mathbf{b}^i/\mathbf{c}^i] \in \mathbb{C}^*$ denote the determinant of the transition matrix between those bases (the entries of this matrix are coordinates of vectors in $d_{i+1}(\mathbf{b}^{i+1})\widetilde{\mathbf{h}}^i\mathbf{b}^i$ with respect to $\mathbf{c}^i$). The \emph{sign-determined Reidemeister torsion} of $C_*$ (with respect to the bases $\mathbf{c}^*$ and $\mathbf{h}^*$) is the following alternating product (see~\cite[Definition 3.1]{Turaev:2000}):
\begin{equation}
\label{Def:RTorsion}
\mathrm{Tor}(C_*, \mathbf{c}^*, \mathbf{h}^*) = (-1)^{|C_*|} \cdot  \prod_{i=0}^n [d_{i+1}(\mathbf{b}^{i+1})\widetilde{\mathbf{h}}^i\mathbf{b}^i/\mathbf{c}^i]^{(-1)^{i+1}} \in \mathbb{C}^*.
\end{equation}
Here  $$|C_*| = \sum_{k\geqslant 0} \alpha_k(C_*) \beta_k(C_*),$$ where $\alpha_i(C_*) = \sum_{k=0}^i \dim C_k$ and  $\beta_i(C_*)  = \sum_{k=0}^i \dim H_k$.

The torsion $\mathrm{Tor}(C_*, \mathbf{c}^*, \mathbf{h}^*)$ does not depend on the choices of $\mathbf{b}^i$ and $\widetilde{\mathbf{h}}^i$. 
Note that if $C_*$ is acyclic (i.e. if $H_i = 0$ for all $i$), then $|C_*| = 0$.

%It does only depend on the equivalence classes of $\mathbf{c}^i$ and $\mathbf{h}^i$. More precisely, if $\mathbf{c'}^i$ is another basis of $C_i$ and $\mathbf{h'}^i$ another one of $H_i$, then we have the so-called \emph{basis change formula}
%\begin{equation}
%\label{EQ:changementdebase}
%\frac{\mathrm{tor}(C_*, \mathbf{c'}^*, \mathbf{h'}^*)}{\mathrm{tor}(C_*, \mathbf{c}^*, \mathbf{h}^*)} = \prod_{i=0}^n \left( \frac{[{\mathbf{c}'}^i/\mathbf{c}^i]}{[{\mathbf{h}'}^i/\mathbf{h}^i]}\right)^{(-1)^i}.
%\end{equation}

\subsubsection*{Torsion of a CW-complex}
Let $W$ be a finite CW-complex and $\rho$ be an $\SL$-re\-pre\-sen\-ta\-tion of $\pi_1(W)$. We define the $\sll_{\rho}$-twisted chain complex of $W$ to be
\[
C_*(W; \sll_\rho) = C_*(\widetilde{W}; \ZZ) \otimes_{\ZZ[\pi_1(W)]} \sll_\rho.
\]
Here $C_*(\widetilde{W}; \ZZ)$ is the complex of the universal cover with integer coefficients which is in fact a $\ZZ[\pi_1(W)]$-module (via the action of $\pi_1(W)$ on $\widetilde{W}$ as the covering group), and $\sll_\rho$ denotes the $\ZZ[\pi_1(W)]$-module via the composition $Ad \circ \rho$, where $Ad \colon \SL \to \mathrm{Aut}(\sll), A \mapsto Ad_A$, is the adjoint representation. The chain complex $C_*(W; \sll_\rho)$ computes the {$\sll_\rho$-twisted homology} of $W$ which we denote as $H_*^\rho(W) = H_i(W; Ad \circ \rho)$.

Let $\left\{e^{(i)}_1, \ldots, e^{(i)}_{n_i}\right\}$ be the set of $i$-dimensional cells of $W$. We lift them to the universal cover and we choose an arbitrary order and an arbitrary orientation for the cells $\left\{ {\tilde{e}^{(i)}_1, \ldots, \tilde{e}^{(i)}_{n_i}} \right\}$. If $\mathcal{B} = \{\mathbf{a}, \mathbf{b}, \mathbf{c}\}$ is an orthonormal basis of $\sll$, then we consider the corresponding  basis over $\CC$
$$\mathbf{c}^{i}_{\mathcal{B}} = \left\{ \tilde{e}^{(i)}_{1} \otimes \mathbf{a}, \tilde{e}^{(i)}_{1} \otimes \mathbf{b}, \tilde{e}^{(i)}_{1} \otimes \mathbf{c}, \ldots, \tilde{e}^{(i)}_{n_i}\otimes \mathbf{a}, \tilde{e}^{(i)}_{n_i} \otimes \mathbf{b}, \tilde{e}^{(i)}_{n_i}\otimes \mathbf{c}\right\}$$ of $C_i(W; \sll_\rho) = C_*(\widetilde{W}; \ZZ) \otimes_{\ZZ[\pi_1(W)]} \sll_\rho$. Now choosing for each $i$ a basis $\mathbf{h}^{i}$ of the $\sll_\rho$-twisted homology $H_i^\rho(W)$, we can compute the torsion
$$\mathrm{Tor}(C_*(W; \sll_\rho), \mathbf{c}^*_{\mathcal{B}}, \mathbf{h}^{*}) \in \CC^*.$$

The cells $\left\{ \tilde{e}^{(i)}_j \right\}^{}_{0 \leqslant i \leqslant \dim W, 1 \leqslant j \leqslant n_i}$ are in one--to--one correspondence with the cells of $W$, their order and orientation induce an order and an orientation for the cells $\left\{ e^{(i)}_j \right\}^{}_{0 \leqslant i \leqslant \dim W, 1 \leqslant j \leqslant n_i}$. Again, corresponding to these choices, we get a basis $c^i$ over $\IR$ of $C_i(W; \IR)$. 

Choose an \emph{homology orientation} of $W$, which is an orientation of the real vector space $H_*(W; \IR) = \bigoplus_{i\geqslant 0} H_i(W; \IR)$. Let $\mathfrak{o}$ denote this chosen orientation. Provide each vector space $H_i(W; \IR)$ with a reference basis $h^i$ such that the basis $\left\{ {h^0, \ldots, h^{\dim W}} \right\}$ of $H_*(W; \IR)$ is {positively oriented} with respect to $\mathfrak{o}$. Compute the sign--determined Reidemeister torsion $\mathrm{Tor}(C_*(W; \IR), c^*, h^{*}) \in \IR^*$ of the resulting based and homology based chain complex and consider its sign $$\tau_0 = \mathrm{sgn}\left(\mathrm{Tor}(C_*(W; \IR), c^*, h^{*})\right) \in \{\pm 1\}.$$  

We define the (sign--refined)  twisted  Reidemeister torsion of $W$ {(with respect to $\mathbf{h}^{*}$ and $\mathfrak{o}$)} to be
\begin{equation}\label{EQ:TorsionRaff}
\mathrm{TOR}(W; \sll_\rho, \mathbf{h}^{*}, \mathfrak{o}) = \tau_0 \cdot \mathrm{Tor}(C_*(W; \sll_\rho), \mathbf{c}^*_{\mathcal{B}}, \mathbf{h}^{*}) \in \mathbb{C}^*.
\end{equation}
This definition only depends on the combinatorial class of $W$, the conjugacy class of $\rho$, the choice of $\mathbf{h}^{*}$ and the {homology} orientation $\mathfrak{o}$. It is independent of the orthonormal basis $\mathcal{B}$ of $\sll$, of the choice of the lifts $\tilde{e}^{(i)}_j$, and of the choice of the positively oriented basis of $H_*(W; \IR)$. Moreover, it is independent of the order and orientation of the cells (because they appear twice). 

One can prove that $\mathrm{TOR}$ is invariant under cellular subdivision, homeomorphism and simple homotopy equivalences. In fact, it is precisely the sign $(-1)^{|C_*|}$ in Eq.~(\ref{Def:RTorsion}) which ensures all these important invariance properties to hold.

\subsection{Regularity for representations}

In this subsection, we briefly review two notions of regularity {(see~\cite{Heu:2003},~\cite{Porti:1997} and~\cite{JDFibre})}. In the sequel $K \subset S^3$ denotes an oriented knot. 

 %The meridian $\mu$ of $K$ is supposed to be oriented according to the rule $\lk(K, \mu) = +1$, while the longitude $\lambda$ is oriented according to the condition $\mathrm{int}(\mu, \lambda) = +1$. Here $\mathrm{int}(\cdot, \cdot)$ denotes the intersection form on $\bord E_K$. 

%An irreducible representation $\rho \colon \Pi(K) \to \SL$ is called \emph{regular} if $\dim H^1_\rho(E_K) = 1$. 

Observe that for any representation $\rho$, $\dim H_1^\rho(E_K)$ is always greater or equal to $1$. We say that $\rho \in R^\mathrm{irr}(\Pi(K); \SL)$ is \emph{regular} if $\dim H_1^\rho(E_K) = 1$. This notion is invariant by conjugation and thus it is well defined for irreducible characters.

\begin{example}
For the trefoil knot and for the figure eight knot, one can prove that each irreducible representation of its group in $\SL$ is regular (see~\cite{JDFourier} and~\cite{Porti:1997}). 
\end{example}

Note that for a regular representation $\rho\colon \Pi(K) \to \SL$, we have $$\dim H_1^\rho(E_K) = 1, \; \dim H_2^\rho(E_K) = 1 \text{ and } H_j^\rho(E_K) = 0 \text{ for all } j \ne 1, 2.$$ 

Let $\gamma$ be a simple closed unoriented curve in $\bord E_K$. Among irreducible representations we focus on the $\gamma$-regular ones. We say that an irreducible representation $\rho : \Pi(K) \to \SL$ is \emph{$\gamma$-regular}, if (see~\cite[Definition 3.21]{Porti:1997}):
\begin{enumerate}
  \item the inclusion $\iota \colon \gamma \hookrightarrow E_K$ induces a \emph{surjective} map $$\iota_* \colon H^\rho_1(\gamma) \to H^\rho_1(E_K),$$
  \item if $\tr(\rho(\pi_1(\bord E_K))) \subset \{\pm 2\}$, then $\rho(\gamma) \ne \pm \I$.
\end{enumerate} 
It is easy to see that this notion is invariant by conjugation, thus the notion of $\gamma$-regularity is well-defined for irreducible characters. Also observe that a $\gamma$-regular representation is necessarily regular (the converse is false in general for an arbitrary curve).

\begin{example}\label{ExReg}
For the trefoil knot, all irreducible representations of its group in $\SL$ are $\lambda$-regular (see~\cite{JDFourier}).

For the figure eight knot, one can prove that each irreducible representation of its group in $\SL$ is $\lambda$-regular \emph{except} two.
%{Question: Is it the case for all twist knots\,? I think (hope) that it is also the case... We must write a proposition on that !}
%\textcolor{blue}{Yoshikazu wrote:  This statement do not hold for all twist knots. We need some corrections.}
\end{example}

We close this section with an important fact concerning hyperbolic knots.
\begin{fact}[\cite{Porti:1997}]\label{holonomy}
Let $K$ be a hyperbolic knot and 
consider the holonomy representation $\rho_0$ associated to the complete hyperbolic structure. 
Let $\gamma$ be any simple closed curve in the boundary of $E_K$ such that $\rho_0(\gamma) \ne \pm \I$, then $\rho_0$ is $\gamma$-regular.
\end{fact}

In particular, for a hyperbolic knot the holonomy representation $\rho_0$ is always $\mu$-regular and $\lambda$-regular.

Applying  \cite[Proposition 3.26]{Porti:1997} to a hyperbolic knot exterior $E_K$,  
we obtain that for any simple closed curve $\gamma$, irreducible and non-$\gamma$-regular characters are contained in the set of zeros of the differential of the trace--function $I_{\gamma}$. 

\begin{remark}
Since the trace--function $I_{\gamma}$ is a regular function 
on the character variety,
the set of irreducible and non-$\gamma$-regular characters is discrete 
on the components where $I_{\gamma}$ is nonconstant.
\end{remark}

If $K$ is a hyperbolic knot, then the character of a complete holonomy representation is contained in a $1$-dimensional irreducible component $X_0(\Pi(K))$ of $X(\Pi(K))$, which satisfies the following condition: if a simple closed curve $\gamma$ in $\partial E_K$ represents any nontrivial element of $\Pi(K)$ then the trace--function $I_{\gamma}$ is nonconstant on $X_0(\Pi(K))$ (see \cite[Corollary 4.5.2]{HandBookGT}).
In particular, irreducible characters near the character of a complete holonomy representation
are $\mu$-regular and $\lambda$-regular.

%Here is an alternative formulation, see~\cite[Proposition 3]{JDFibre}. Fix a generator $P^\rho$ of $H^0_\rho(\bord E_K)$. The inclusion $\iota \colon \gamma \hookrightarrow E_K$ and the cup product  induce the linear form $f^\rho_{\gamma} \colon H^1_\rho(E_K) \to \mathbb{C}$. We explicitly have $$f^\rho_{{\gamma}}(v) = B_{\sll} \left( {P^\rho, v(\gamma)}\right), \text{ for all } v \in  H^1_\rho(E_K).$$ 

%\begin{proposition}[Proposition 3 of \cite{JDFibre}]
%An irreducible representation $\rho$ of $\Pi(K)$ in $\SL$ is $\gamma$-regular if and only if the linear form $f^\rho_{\gamma} \colon H^1_\rho(E_K) \to \mathbb{C}$ is an isomorphism.
%\end{proposition}

\subsection{Review on the non--abelian Reidemeister torsion for knot exteriors}
\label{Torsion}
This subsection gives a detailed review of the constructions made in~\cite[Section~6]{JDFourier}. In particular, we shall explain how to construct distinguished bases for the twisted homology groups of knot exteriors. 

\subsubsection*{Canonical homology orientation of knot exteriors}
We equip the exterior of $K$ with its \emph{canonical homology orientation} defined as follows (see~\cite[Section V.3]{Turaev:2002}). We have 
$$H_*(E_K; \IR) = H_0(E_K; \IR) \oplus H_1(E_K; \IR)$$ 
and we base this $\IR$-vector space with $\{ \lbrack \! \lbrack pt \rbrack \! \rbrack, \lbrack \! \lbrack \mu \rbrack \! \rbrack\}$. Here $\lbrack \! \lbrack pt \rbrack \! \rbrack$ is the homology class of a point, and $\lbrack \! \lbrack \mu \rbrack \! \rbrack$ is the homology class of the meridian $\mu$ of $K$. This reference basis of $H_*(E_K; \IR)$ induces the so--called canonical homology orientation of $E_K$. In the sequel, we let $\mathfrak{o}$ denote the canonical homology orientation of $E_K$.

\subsubsection*{How to construct natural bases of the twisted homology}

Let $\rho$ be a regular $\SL$-representation of $\Pi(K)$ and {\emph{fix}} a generator $P^\rho$ of $H_0^\rho(\bord E_K)$ (i.e. $P^\rho$ is an element in $\sll$ such that $Ad_{\rho(g)}(P^\rho) = P^\rho$ for all $g \in \pi_1(\bord E_K)$).  

The canonical inclusion $i\colon \bord E_K \to E_K$ induces {(see~\cite[Corollary 3.23]{Porti:1997})} an isomorphism $i_*\colon H_2^\rho(\bord E_K) \to H_2^\rho(E_K)$. Moreover, one can prove that {(see~\cite[Proposition 3.18]{Porti:1997})}
$$H_2^\rho(\bord E_K) \cong H_2(\bord E_K; \ZZ) \otimes \CC.$$
 More precisely, let $\lbrack \! \lbrack \bord E_K \rbrack \! \rbrack \in H_2(\bord E_K; \ZZ)$ be the fundamental class induced by the orientation of $\bord E_K$, one has $H_2^\rho(\bord E_K) = \CC\left[\lbrack \! \lbrack \bord E_K \rbrack \! \rbrack \otimes P^\rho\right]$.

 The \emph{reference generator} of $H_2^\rho(E_K)$ is defined by 
\begin{equation}\label{EQ:Defh2}
h_{(2)}^\rho 
= i_*([\lbrack \! \lbrack \bord E_K \rbrack \! \rbrack \otimes P^\rho]).
\end{equation}

Let $\rho$ be a $\lambda$-regular representation of $\Pi(K)$. The \emph{reference generator} of the first twisted homology group $H_1^\rho(E_K)$ is defined by
\begin{equation}\label{EQ:Defh1}
h_{(1)}^\rho(\lambda) 
= \iota_*\left(\left[\lbrack \! \lbrack \lambda \rbrack \! \rbrack  \otimes P^\rho\right]\right).
\end{equation}

\begin{remark}
The generator $h_{(1)}^\rho(\lambda)$ of  $H_1^\rho(E_K)$ depends on the orientation of $\lambda$. If we change the orientation of  the longitude $\lambda$ in Eq.~(\ref{EQ:Defh1}),  then the generator is change into its reverse.
\end{remark}

\begin{remark}
Note that $H_i^\rho(E_K)$ is isomorphic to the dual space of the $\sll_\rho$-twisted cohomology $H^i_\rho(E_K) = H^i(E_K; Ad \circ \rho)$. Reference elements defined in Eqs.~(\ref{EQ:Defh2}) and (\ref{EQ:Defh1}) are dual from the ones defined in~\cite[Section~3.4]{JDFibre}.
\end{remark}

\subsubsection*{The Reidemeister torsion for knot exteriors}

Let $\rho \colon \Pi(K) \to \SL$ be a $\lambda$-regular representation. The \emph{Reidemeister torsion $\mathbb{T}^K_\lambda$} at $\rho$ is  
defined to be 
\begin{equation}\label{Tordef}
\mathbb{T}^K_\lambda(\rho)  = \mathrm{TOR}\left( {E_K; \sll_\rho, \{h_{(1)}^\rho(\lambda), h_{(2)}^\rho\}, \mathfrak{o}} \right) \in \CC^*.
\end{equation}
It is an invariant of knots. Moreover, if $\rho_1$ and $\rho_2$ are two $\lambda$-regular representations which have the same character then $\mathbb{T}^K_\lambda(\rho_1) = \mathbb{T}^K_\lambda(\rho_2)$. Thus, $\mathbb{T}^K_\lambda$ defines a map on the set $X^{\mathrm{irr}}_\lambda(\Pi(K)) = \{\chi \in X^{\mathrm{irr}}(\Pi(K)) \; |\; \chi \text{ is } \lambda\text{-regular}\}$ {of $\lambda$-regular characters}.

\begin{remark}\label{rmkinv}
The Reidemeister torsion $\mathbb{T}^K_\lambda(\rho)$ defined in Eq.~(\ref{Tordef}) is exactly the inverse of the one considered in~\cite{JDFibre}.
\end{remark}

{In the following remark we discuss sign properties of the torsion $\mathbb{T}^K_\lambda$.}
\begin{remark}\label{rmsgn}
\begin{enumerate}
  \item {The torsion $\mathbb{T}^{K}_{\lambda}$ does not depend on the orientation of $K$ (see~\cite[Proposition 3.4]{JDFibre}).}
  \item {Let $K^*$ denote the mirror image of $K$, and $\lambda^*$ be its preferred longitude as defined in Section~\ref{section:twist_knots}. The character varieties of $E_K$ and $E_{K^*}$ are same and
$\mathbb{T}^{K^*}_{\lambda^*} =  \mathbb{T}^K_\lambda.$}

{Indeed, if we take the mirror image of $K$, then the orientation of the ambient $3$-sphere is reversed, so the orientation of $E_{K^*}$ is the opposite of the one of $E_K$. As a consequence the generator of $H^\rho_2(E_{K^*})$ defined in Eq.~(\ref{EQ:Defh2}) is the opposite of the one of $H^\rho_2(E_{K})$. On the other hand, the meridian $\mu^*$ of $K^*$ is  the inverse of the one of $K$ whereas longitudes are same. So, the homology orientation of $K^*$ is the opposite of the one of $K$ whereas generators of the twisted $H_1$ are same. As a consequence $\mathbb{T}^{K^*}_{\lambda^*} =  \mathbb{T}^K_\lambda$.}
\end{enumerate}

\end{remark}

\subsection{Review on the non--abelian Reidemeister torsion polynomial}
\label{TorsionPoly}
To compute the non--abelian Reidemeister torsion for twist knots, we use techniques developed by the third author in~\cite{YY2}. In fact, we compute a more general invariant of knots called the non--abelian Reidemeister torsion polynomial. It is a sort of {twisted} Alexander polynomial invariant (but with non--abelian twisted coefficients) whose \lq\lq derivative coefficient" at $t = 1$ is exactly $\mathbb{T}^K_\lambda$.

\subsubsection*{Definitions} Let $W$ be a finite CW--complex.
We regard $\Z$ as a multiplicative group which is generated by one variable $t$.
Let $\alpha$ be a surjective homomorphism from $\pi_1(W)$ to $\Z = \langle t \rangle$. 

If $\rho$ is an $\SL$-representation of $\pi_1(W)$, we define the $\tsll_\rho$-twisted chain complex of $W$ to be
\[
    C_*(W;\tsll_\rho) = C_*(\widetilde W;\Z) 
                           \otimes_{Ad \circ \rho \otimes \alpha} 
                           \left(\sll \otimes \C(t)\right),
\]
where $\sigma \cdot \gamma \otimes v \otimes f$ is identified with 
$\sigma \otimes Ad_{\rho(\gamma)}(v) \otimes f \cdot t^{\alpha(\gamma)}$.

The sign--refined Reidemeister torsion of $W$ with respect to this $\tsll_\rho$-twisted {coefficients} is defined to be (compare with Eq.~(\ref{EQ:TorsionRaff}))
\[
\mathrm{TOR}(W;\tsll_\rho, \mathbf{h}^*, \mathfrak{o})
=
\tau_0 \cdot \mathrm{Tor}(C_*(W;\tsll_\rho), \mathbf{c}^*_{\mathcal{B}}, \mathbf{h}^*) \in \C(t)^{*}.
\]
Note that
$\mathrm{TOR}(W;\tsll_\rho, \mathbf{h}^*, \mathfrak{o})$ is --- as the Alexander polynomial --- determined up to 
a factor $t^m$ where $m \in \Z$.

Next we turn back to knots exteriors. From now on, we suppose that the CW--complex $W$ is $E_K$ and that the homomorphism $\alpha : \Pi(K) \to \Z$ is the abelianization. 
From \cite[Proposition 3.1.1]{YY1}, we know that if $\rho$ is $\lambda$-regular, then all homology groups $H_*(E_K;\tsll_\rho)$ vanishes.
So if $\rho$ is $\lambda$-regular, then we define the non--abelian Reidemeister torsion polynomial at $\rho$ to be
\begin{equation}\label{TORpol}
\T^{K}_{\lambda}(\rho)
=
\mathrm{TOR}(W;\tsll_\rho, \emptyset, \mathfrak{o}) \in \C(t)^*.  
\end{equation}
The torsion in Eq.~(\ref{TORpol}) is also determined up to a factor $t^m$ where $m \in \Z$.
It is also shown in \cite[Theorem 3.1.2]{YY1} that 
\[
\mathbb{T}^K_{\lambda}(\rho) = - \lim_{t \to 1}
\frac{\T^K_{\lambda}(\rho)}{t-1}.
\]

\begin{remark}
It is shown by T. Kitano \cite[Theorem A]{Kitano} that $\T^{K}_{\lambda}(\rho)$ agree with the twisted Alexander invariant 
of $K$ and {$Ad \circ \rho \otimes \alpha$}. 
\end{remark}

\subsubsection*{How to compute $\T^K_{\lambda}(\rho)$ from Fox--calculus} 
Here we review a description of $\T^K_{\lambda}(\rho)$ from a Wirtinger presentation of $\Pi(K)$. This description comes from some results by T. Kitano~\cite{Kitano}.
For simplicity, write $\Phi$ for $(Ad \circ \rho) \otimes \alpha$.
Choose a Wirtinger presentation
\begin{equation}\label{presentation}
\Pi(K)
=
\langle x_1, \ldots, x_k \,|\, r_1,\ldots, r_{k-1} \rangle
\end{equation}
of $\Pi(K)$.
Let $W_K$ be the $2$-dimensional CW--complex constructed from the presentation~(\ref{presentation}) in the usual way. The $0$-skeleton of $W_K$ consists of a single $0$-cell $pt$, the $1$-skeleton is a wedge of
$k$ oriented $1$-cells $x_1, \ldots x_k$ and the $2$-skeleton consists of $(k-1)$ $2$-cells $D_1, \ldots, D_{k-1}$ with
attaching maps given by the relations $r_1, \ldots, r_{k-1}$ of presentation~(\ref{presentation}).

 F. Waldhausen proved~\cite{W} that the Whitehead 
group of a knot group is trivial. As a result, $W_K$ has the same simple 
homotopy type as $E_K$. So, the CW--complex $W_K$ can be used to compute  the non--abelian Reidemeister torsion polynomial defined in Eq.~(\ref{TORpol}). 
Therefore it is enough to consider the Reidemeister torsion of the $\tsll_\rho$-twisted chain complex $C_*\left(W_K;\tsll_\rho\right)$.

The twisted complex $C_*(W_K;\tsll_\rho)$ thus becomes:
\begin{equation}\label{twistcomplex}
\xymatrix@1@-.5pc{0 \ar[r] &  {\left({\sll\otimes \C(t)}\right)}_{}^{k-1} 
\ar[r]^-{\partial_2} &  {\left({\sll \otimes \C(t)}\right)}_{}^k 
\ar[r]^-{\partial_1}& \sll \otimes \C(t) \ar[r] & 0.}
\end{equation}
Here we briefly denote the $l$-times direct sum of $\sll \otimes \C(t)$ by $(\sll \otimes \C(t))^l$. In complex~(\ref{twistcomplex}), we have 
\[
\partial_1 
=
\left( {\Phi(x_1 - 1), \Phi(x_2 - 1), \ldots, \Phi(x_k - 1)} \right).
\]
and $\partial_2$ is expressed using the Fox differential calculus and the action given by $\Phi = (Ad \circ \rho) \otimes \alpha$:
\begin{equation}\label{Eq:MatrixD2}
\partial_2 =
\left(
\begin{array}{ccc}
\Phi(\frac{\partial r_1}{\partial x_1}) & \ldots & \Phi(\frac{\partial r_{k-1}}{\partial x_1}) \\
\vdots & \ddots & \vdots \\
\Phi(\frac{\partial r_1}{\partial x_k}) & \ldots & \Phi(\frac{\partial r_{k-1}}{\partial x_k}) 
\end{array}
\right)
\end{equation}

Let $A^i_{K, Ad\circ \rho}$ denote the $3(k-1)\times 3(k-1)$--matrix obtained from the matrix in Eq.~(\ref{Eq:MatrixD2}) by deleting its $i$-th row. The torsion polynomial $\T^K_{\lambda}(\rho)$ defined in Eq.~(\ref{TORpol}) can be described, up to a factor $t^m$ $(m\in \Z)$, as follows  (for more details see \cite{KL, Kitano}):
\begin{equation}\label{torsionK}
\T^K_{\lambda}(\rho) = \tau_0 \cdot
\frac{\det A^i_{K, Ad\circ \rho}}{\det(\Phi(x_i-1))}.
\end{equation}
This rational function has the first order zero at $t=1$ \cite[Theorem 3.1.2]{YY1}.
The non--abelian Reidemeister torsion ${\mathbb T}^K_{\lambda}(\rho)$
is expressed as 
\begin{equation}\label{eqn:torsion_eqn}
{\mathbb T}^K_{\lambda}(\rho)
= -\lim_{t \to 1}
\frac{\T^K_{\lambda}(\rho)}{t-1}
=
-\lim_{t \to 1}
\left(
\tau_0 \cdot
\frac{\det A^i_{K, Ad\circ \rho}}{(t-1)\det(\Phi(x_i-1))}
\right).
\end{equation}

\begin{remark}
From~\cite[Proposition 4.3.1]{YY1},
we can see that 
the non--abelian Reidemeister torsion ${\mathbb T}^K_{\lambda}$ associated to 
a two--bridge knot $K$ is a rational function in $s+1/s$ and $u$, where $(s, u)$ is a solution of Riley's equation $\phi_K(s,u) = 0$.
In particular, if we consider the case for $s=1$, then 
the Reidemeister torsion ${\mathbb T}^K_{\lambda}$ is a rational function of $u$.
The variable $u$ satisfies Riley's equation $\phi_K(1, u) = 0$.
For a hyperbolic twist knot $K$, $u$ is expressed in terms of the cusp shape.
Thus the non--abelian Reidemeister torsion ${\mathbb T}^K_{\lambda}$ at the holonomy $\rho_0$ is also a rational function 
in the cusp shape of $K$.
\end{remark}

%%%%%%%%%%%%%%%%%%%%%%%%%%%%%%%%%%%
\section{The non--abelian Reidemeister torsion for twist knots}
\label{section:results}
%%%%%%%%%%%%%%%%%%%%%%%%%%%%%%%%%%%
In this section, we compute the non--abelian Reidemeister torsion for twist knots. 
Since there exists an isomorphism between the knot groups $\Pi(J(2, 2m+1))$ and $\Pi(J(2, -2m))$ (see Remark~\ref{remarktwist}), it is enough for us to make the computations in the case of even twist knots $K = J(2, 2m)$, $m \in \ZZ$.
The method used is the following. We will make the computation at the acyclic level, i.e. compute the torsion polynomial $\T^K_{\lambda}(\rho)$, and next apply~\cite[Theorem 3.1.2]{YY1} to obtain $\mathbb{T}^K_\lambda(\rho)$ (see Eq.~(\ref{eqn:torsion_eqn})). 

\begin{remark}
 Remark~\ref{remarktwist}, item (5), says that the knot $J(2, -2m)$ is obtained by surgery on the Whitehead link. Of course the surgery formula for the Reidemeister torsion (see e.g.~\cite[Theorem 4.1 (iii)]{Porti:1997}) theoretically gives a formula for the non--abelian torsion for $J(2, -2m)$, but unfortunately it is difficult to extract from it an explicit formula as we are interested in this paper.
\end{remark}

\subsection{The non--abelian Reidemeister torsion for  even twist knots}
We calculate the non--abelian Reidemeister torsion for  even twist knots $J(2, 2m)$ where $m$ is an integer.

\subsubsection{Preliminaries}
Following Section \ref{ReviewCharK} and using Riley's method, we can parametrize a non--abelian $\SL$-representation $\rho$ by two parameters $u$ and $s$  as follows:
\[
\rho(x) = 
\left(
\begin{array}{cc}
\sqrt{s} & 1/\sqrt{s} \\
0 & 1/\sqrt{s}
\end{array}
\right),\,
\rho(y) = 
\left(
\begin{array}{cc}
\sqrt{s} & 0 \\
-\sqrt{s} u & 1/\sqrt{s}
\end{array}
\right),
\]
where $s$ and $u$ satisfy Riley's equation $\phi_{J(2, 2m)}(s, u) = 0$. Besides, the Riley polynomial for twist knots is such that: 
\begin{equation}\label{eqn:Riley_poly}
\phi_{J(2, 2m)} (s,u)
=
\frac{(s+s^{-1}-1-u)(\xi_{+}^m - \xi_{-}^m) -(\xi_{+}^{m-1} - \xi_{-}^{m-1})}{\xi_{+} - \xi_{-}},
\end{equation}
where $\xi_{\pm}$ are the \emph{eigenvalues}  of the matrix $\rho(w) = \rho([y, x^{-1}])$ given by 
\begin{equation}\label{ev}
\xi_{\pm}   = \frac{1}{2} \left[ {u^2+(2-s-s^{-1})u+2 \pm \sqrt{(u^2+(2-s-s^{-1})u+4)(u^2+(2-s-s^{-1})u)}}\right].
\end{equation}

% \begin{figure}[!htbp]
% \begin{center}
% \scalebox{.4}{\includegraphics{GraphCuspShapes}}
% \end{center}
% \caption{Graph of the cusp shape of $J(2,-2m)$.}
% \label{graph:cusp_shapes}
% \end{figure}

\subsubsection{Statement of the result}
\label{statement}

\begin{notation*}
Let $\alpha_1$, $\alpha_2$, $\beta_1$, $\beta_2$, $c$ and  $t_m$ be as follows:
% \begin{align*}
% c &= c(u,s) = u+1-s^{-1}; \\
% \alpha_1 &= (\xi_{-} - 1)(\xi_{+} + s) + (s-1)^2 u-su^2; \\
% \alpha_2 &=  (1-su-\xi_{+})\left(1+\frac{\xi_{-} - s}{c}\right); \\
% \beta_1   &= (\xi_{+}-1)(\xi_{-} + s) + (s-1)^2 u-su^2; \\
% \beta_2   &=  (1-su-\xi_{-})\left(1+\frac{\xi_{+} - s}{c} \right); \\
% t_m  &= \frac{\xi_{+}^m -\xi_{-}^m}{\xi_{+} -\xi_{-}}.% t_m
% \end{align*}
\begin{itemize}
\item[] $c = c(s,u) = u+1-s^{-1}$; 
\item[] $\alpha_1 = (\xi_{-} - 1)(\xi_{+} + s) + (s-1)^2 u-su^2$; 
\item[] $\alpha_2 =  (1-su-\xi_{+})\left(1+ (\xi_{+} - s)/c\right)$; 
\item[] $\beta_1  = (\xi_{+}-1)(\xi_{-} + s) + (s-1)^2 u-su^2$;
\item[] $\beta_2  =  (1-su-\xi_{-})\left(1+(\xi_{-} - s)/c \right)$;
\item[] $t_m  = (\xi_{+}^m -\xi_{-}^m)/(\xi_{+} -\xi_{-})$.
\end{itemize}
\end{notation*}

\begin{remark}
Using such notation, the Riley polynomial of the twist knot $J(2, 2m)$ becomes:
\[
\phi_{J(2, 2m)} (s,u)
=
(s-c)t_m - t_{m-1}.
\]
\end{remark}

With this notation in mind we can write down the general formula for the non--abelian Reidemeister torsion for twist knots.

\begin{theorem}\label{theorem_J(2,2m)}
Let $m$ be a positive integer.
\begin{enumerate}
  \item  The Reidemeister torsion ${\mathbb T}^{J(2, 2m)}_{\lambda}(\rho)$ satisfies the following formula:
\begin{equation}\label{m_pos}
{\mathbb T}^{J(2, 2m)}_{\lambda}(\rho)
=
\frac{\tau_0}{s+s^{-1} - 2} 
\left[
    C_1(m) \xi_{+}^{m-1} t_m + C_2(m) \xi_{-}^{m-1} t_m + C_3(m)
\right].
\end{equation}
\item  Similarly, we have
\begin{equation}\label{m_neg}
{\mathbb T}^{J(2, -2m)}_{\lambda}(\rho)
=
\frac{\tau_0}{s+s^{-1} - 2} 
\left[
    -C_1(-m) \xi_{+}^{-m-1} t_m - C_2(-m) \xi_{-}^{-m-1} t_m + C_3(-m)
\right].
\end{equation}
In that two formulas we have:
\begin{align*}%%%% C_1
C_1(m)
&= 
\frac{1}{(\xi_{+} - \xi_{-})^2} 
\left\{
\frac{1}{s} \alpha_1^2 (3m+1) + \frac{m}{s}\beta_1^2 (\xi_{+}^2+1) - m(\xi_{+}-\xi_{-})^2(s+\frac{1}{s}+1) 
\right\}\\
& \quad 
    -\frac{m}{(\xi_{+}-\xi_{-})^4}
    \left\{
    \left( c(1-\xi_{+}) + \frac{\alpha_1}{s}\right)^2  
    \left( (\xi_{+}-\xi_{-})^2(s+\frac{1}{s}+1) - \frac{\alpha_1^2+\beta_1^2}{s}\right)
    \right.\\
& \left. \qquad \qquad \qquad
     +\frac{2\alpha_1\beta_2}{s}\left(c(1-\xi_{+})+\frac{\alpha_1}{s}\right)\left(c(1-\xi_{+})+\frac{\alpha_2}{s}\right) \right\} \\
& \quad
    -\frac{m}{(\xi_{+}-\xi_{-})^4}
    \left\{
       \frac{\alpha_1^2}{s}
        \left( 
        (\xi_{+} - \xi_{-})^2(u^2+4u+3) - \left( c(1-\xi_{+}) + \frac{\alpha_1}{s}\right)^2 
        \right.
    \right.\\
&        \qquad \qquad \qquad \qquad \qquad
    \left.
        \left.
        -\left(c(1-\xi_{-}) + \frac{\beta_1}{s}\right)^2
        \right)
        \right. \\
& \qquad \qquad \qquad \left. 
        + \frac{2\alpha_1\alpha_2}{s}
        \left(
        c(1-\xi_{+}) +\frac{\alpha_1}{s}
        \right)
        \left(
        c(1-\xi_{-}) +\frac{\beta_2}{s}
        \right) 
        \right\},
\end{align*}
\begin{align*}%%% C_2
C_2(m)
&=
\frac{1}{(\xi_{+} -\xi_{-})^2}
\left\{
\frac{1}{s} \beta_1^2 (3m+1) + \frac{m}{s}\alpha_1^2 (\xi_{-}^2+1) - m(\xi_{+}-\xi_{-})^2(s+\frac{1}{s}+1) 
\right\} \\
& \quad
-\frac{m}{(\xi_{+}-\xi_{-})^4}
    \left\{
    \left(
    c(1-\xi_{-})+\frac{\beta_1}{s}
    \right)^2
    \left(
    (\xi_{+}-\xi_{-})^2(s+\frac{1}{s}+1)-\frac{\alpha_1^2 + \beta_1^2}{s}
    \right)
    \right.\\
& \qquad \qquad \qquad \left. 
    +\frac{2\alpha_2 \beta_1}{s}
    \left(
    c(1-\xi_{-})+\frac{\beta_1}{s}
    \right)
    \left(
    c(1-\xi_{-})+\frac{\beta_2}{s}
    \right)
    \right\}\\
& \quad
-\frac{m}{(\xi_{+}-\xi_{-})^4}
    \left\{
    \frac{\beta_1^2}{s} 
      \left( 
      (\xi_{+} - \xi_{-})^2(u^2+4u+3) - \left(c(1-\xi_{+}) + \frac{\alpha_1 }{s}\right)^2 
      \right.
    \right.\\
& \qquad \qquad \qquad \qquad \qquad
    \left.
      \left.
         - \left(c(1-\xi_{-}) + \frac{\beta_1}{s}\right)^2 
      \right) 
    \right. \\
& \qquad \qquad \qquad 
        \left.
        +\frac{2\beta_1\beta_2}{s}
        \left(c(1-\xi_{-})+\frac{\beta_1}{s}\right)\left(c(1-\xi_{+})+\frac{\alpha_2}{s}\right)
        \right\},
\end{align*}
\begin{align*}%%%%%% C_3
C_3(m)
&=
\frac{m}{(\xi_{+} -\xi_{-})^2}
\left\{
(\xi_{+}-\xi_{-})^2(s+\frac{1}{s}+1)-\frac{\alpha_1^2 + \beta_1^2}{s}
\right\} \\
&\quad+
\frac{t_m^2}{(\xi_{+} -\xi_{-})^2}
\left\{
4(\xi_{+}-\xi_{-})^2(s+\frac{1}{s}+1)-\frac{5(\alpha_1^2+\beta_1^2)}{s}
\right\} \\
&\quad -
 \frac{t_m^2}{(\xi_{+}-\xi_{-})^4}
\left\{
    \frac{1}{s}
    \left( 
    c(1-\xi_{+})\beta_1+\frac{\alpha_1 \beta_1}{s}
    \right)^2
    -\frac{1}{s}
    \left(
    c(1-\xi_{+})\beta_2+\frac{\alpha_2 \beta_2}{s}
    \right)^2
\right.\\
& \left.\qquad \qquad \qquad
    +\frac{1}{s}
    \left(
    c(1-\xi_{-})\alpha_1+\frac{\alpha_1 \beta_1}{s}
    \right)^2
    -\frac{1}{s}
    \left(
    c(1-\xi_{-})\alpha_2 +\frac{\alpha_2 \beta_2}{s}
    \right)^2
\right\} \\
&\quad +m(s+\frac{1}{s}-2)^2 t_m^2.
\end{align*}
\end{enumerate}
\end{theorem}

\begin{remark}
One can observe that 
${\mathbb T}^{J(2, 2m)}_{\lambda}$ is symmetric in $\xi_{\pm}$.
Together with the fact that $\xi_{+} \cdot \xi_{-} = 1$, 
we can see that ${\mathbb T}^{J(2, 2m)}_{\lambda}$ is in fact a function of $\xi_{+} + \xi_{-} = u^2 + (2-s-s^{-1})u+2$.
\end{remark}

\subsubsection{Proof of Theorem~\ref{theorem_J(2,2m)}}
We make the detailed proof in the case of  $J(2, 2m)$ for $m>0$.

First, recall that the group of $J(2, 2m)$ admits the following presentation (see Fact~\ref{Fact1}): 
\[
\langle x, y \,|\, w^m x  = y w^m \rangle.
\]
Here $w$ is the word $[y, x^{-1}] = yx^{-1}y^{-1}x$.

Before computations, we give an elementary and useful {linear algebra} lemma about trace of matrices in $M_3(\C)$.

\begin{lemma} \label{lemma_sigma_2}
The two following items hold:
\begin{enumerate}
\item
    Let $A$ be in $M_3(\C)$.  Set
\[
\sigma_1(A) = \trace(A), 
\sigma_2(A) = \frac{1}{2} \left( \trace^2(A)-\trace(A^2)\right) \text{ and } \sigma_3(A) = \det (A).
\] 
    We have 
\[
\det(\I + A) = 1 + \sigma_1(A) + \sigma_2(A) + \sigma_3(A)
\]  
    and if $A \in {\rm GL}(3, \C)$,
$$
\sigma_2(A) = \sigma_1(A^{-1})\cdot \det(A).
$$
\item
    If $A = {(a_{i, j})}_{1 \leq i, j \leq 3}$ and $B = {(b_{i, j})}_{1 \leq i, j \leq 3}$ are two matrices in $M_3(\C)$,
    then we have 
    \begin{align*}
\trace(A)\trace(B) - \trace(AB) 
&=
    \left|
    \begin{array}{cc}
    a_{1, 1} & a_{1, 3} \\
    b_{3, 1} & b_{3, 3}
    \end{array}
    \right|
    +
    \left|
    \begin{array}{cc}
    b_{1, 1} & b_{1, 3} \\
    a_{3, 1} & a_{3, 3}
    \end{array}
    \right|\\
& \quad +
    \left|
    \begin{array}{cc}
    a_{2, 2} & a_{2, 3} \\
    b_{3, 2} & b_{3, 3}
    \end{array}
    \right|
    +
    \left|
    \begin{array}{cc}
    b_{2, 2} & b_{2, 3} \\
    a_{3, 2} & a_{3, 3}
    \end{array}
    \right|\\
&\quad +
    \left|
    \begin{array}{cc}
    a_{1, 1} & a_{1, 2} \\
    b_{2, 1} & b_{2, 2}
    \end{array}
    \right|
    +
    \left|
    \begin{array}{cc}
    b_{1, 1} & b_{1, 2} \\
    a_{2, 1} & a_{2, 2}
    \end{array}
    \right|.
    \end{align*}
\end{enumerate}
\end{lemma}
{\begin{remark}
Observe that $\sigma_2(A)$ is the trace of the matrix of cofactors of $A$.
\end{remark}}

\begin{proof}
\phantom{}\hfill
\begin{enumerate}
  \item The first Eq. of item 1 is well-known. From the Cayley--Hamilton identity, we have
\[
A^3 - \sigma_1(A) A^2 + \sigma_2(A) A - \det(A) \I= 0.
\]
% Multiplying this equation by $-\det(A)^{-1}A^{-3}$, we obtain the second Eq. of our first claim.
Multiplying this equation by $A^{-1}$ then taking traces of both sides, we obtain the second Eq. of our first claim.
  \item Second item follows from direct calculations.
\end{enumerate}
\end{proof}

%%%%%%%%%%%%%%%%%%%%%%%%%%%
\noindent{\it Fox--differential calculus for $2m$-twist knots.}
%%%%%%%%%%%%%%%%%%%%%%%%%%%
Since $J(2, 2m)$ is a two--bridge knot,
the non--abelian Reidemeister torsion polynomial
$\T^K_\lambda(\rho)$ associated to $J(2, 2m)$
is expressed as (see Eq.~(\ref{torsionK})):

\begin{equation}\label{torsionJ}
\T^K_\lambda(\rho) = \tau_0
\frac{\det \Phi(\frac{\partial}{\partial x}w^m x w^{-m} y^{-1} )  }
{\det \Phi(y-1)}
\end{equation}
where $\Phi=Ad\circ \rho \otimes \alpha: \Z [\Pi(J(2, 2m))] \to M_3(\C[t, t^{-1}])$.

The following claim gives us the Fox--differential part in the numerator of Eq.~(\ref{torsionJ}).

\begin{claim}\label{ClaimFox}
For $m > 0$, we have:
\begin{equation}\label{EqClaimFox}
\frac{\partial}{\partial x} \left( w^m x w^{-m} y^{-1} \right)  = w^m\left(1 + (1-x)(1 + w^{-1} + \cdots + w^{-m+1})(x^{-1}-x^{-1}y) \right).
\end{equation}
\end{claim}

\begin{proof}[Proof of Claim~\ref{ClaimFox}]
We have:
\begin{align*}
\frac{\partial}{\partial x} \left( w^m x w^{-m} y^{-1} \right) 
&=
\frac{\partial w^m}{\partial x}
+ w^m
+ w^m x (-w^{-m})\frac{\partial w^m}{\partial x} \nonumber \\
&=
w^m
\left(1+(1-x)w^{-m}\frac{\partial w^m}{\partial x}\right). 
\end{align*}
It is easy to see that
\[
\frac{\partial w}{\partial x} = \frac{\partial}{\partial x}\left(yx^{-1}y^{-1}x\right) = yx^{-1}y^{-1} - yx^{-1}.
\] 
Thus
\begin{align*}
w^{-m}\frac{\partial w^m}{\partial x}
&=
(1+ w^{-1} +\cdots + w^{-m+1})w^{-1}\frac{\partial w}{\partial x}
\nonumber\\
&=
(1 + w^{-1} + \cdots + w^{-m+1})(x^{-1}-x^{-1}y),
\end{align*}
which gives us Eq.~(\ref{EqClaimFox}).
\end{proof}

Let $\{E, H, F\}$ be the following usual $\CC$-basis of the Lie algebra $\sll$: 
\[
E=
\left(
\begin{array}{cc}
0 & 1 \\
0 & 0
\end{array}
\right),
H=
\left(
\begin{array}{cc}
1 & 0 \\
0 & -1
\end{array}
\right),
F=
\left(
\begin{array}{cc}
0 & 0 \\
1 & 0
\end{array}
\right).
\]
It is easy to see that the adjoint actions of $x$ and $y$ in the basis $\{ E, H, F \}$ of $\sll$ are given by the following matrices:
\[
X = Ad_{\rho(x)} =
\left(
\begin{array}{ccc}
 s & -2 & -s^{-1} \\
 0 & 1 & s^{-1} \\
 0 & 0 & s^{-1}
\end{array}
\right),\quad
Y = Ad_{\rho(y)} =
\left(
\begin{array}{ccc}
 s & 0 & 0 \\
 su & 1 & 0 \\
 -su^2 & -2u & s^{-1}
\end{array}
\right).
\]

If $W = Ad_{\rho(w)}$, then $\Phi(\frac{\partial}{\partial x}w^m  x w^{-m} y^{-1})$ is given by (see Claim~\ref{ClaimFox}):
\[
W^m 
\left(\I + (\I - tX)(\I+ W^{-1} +\cdots +W^{-m+1})(t^{-1}X^{-1}-X^{-1}Y)\right).
\]

Set $S_m(A) = \I  + A + \cdots + A^{m-1}$, for $A \in \SL$, we finally obtain: 
\begin{equation}\label{PHI}
\Phi(\frac{\partial}{\partial x}w^m  x w^{-m} y^{-1}) = W^m
\left(
\I +(\I -tX)S_m(W^{-1})(t^{-1}X^{-1} - X^{-1}Y)
\right).
\end{equation}

%%%%%%%%%%%%%%%%%%%%%%%%%%%%%%%%%
\noindent {\it Observation about the \lq\lq second differential" of a determinant.}
%%%%%%%%%%%%%%%%%%%%%%%%%%%%%%%%%
We can compute 
the non--abelian Reidemeister torsion for $J(2, 2m)$ 
combining Eqs.~(\ref{eqn:torsion_eqn}) \&~(\ref{PHI}) as follows: 
\begin{equation*}
{\mathbb T}^{J(2, 2m)}_{\lambda}(\rho)
=
-\tau_0
\lim_{t \to 1}
\frac{\det ( \Phi(\frac{\partial}{\partial x}w^m  x w^{-m} y^{-1})) }
            {(t-1)\det ( \Phi(y-1))}.
\end{equation*}

Using the fact that $\det W =1$, Eq.~(\ref{PHI}) gives:
\begin{equation}
\det \Phi(\frac{\partial}{\partial x}w^m  x w^{-m} y^{-1})
=
\det \left(\I +(\I -tX)S_m(W^{-1})(t^{-1} X^{-1} - X^{-1}Y) \right).
\label{content_Z_m}
\end{equation}
If we write $\det(\I + Z_m)$ for
the right hand side of Eq.~(\ref{content_Z_m}),
then 
\begin{equation}\label{Tlim}
{\mathbb T}^{J(2, 2m)}_{\lambda}(\rho) = 
-\tau_0
\lim_{t \to 1}
\frac{\det(\I+Z_m)}{(t-1)^2(t^2 - (s+s^{-1})t+1)}
\end{equation}
thus
\[{\mathbb T}^{J(2, 2m)}_{\lambda}(\rho)  =
\frac{\tau_0}{s+s^{-1} - 2}
\lim_{t \to 1}
\frac{\det(\I+Z_m)}{(t-1)^2}.
\] 
{Moreover, using the first item of Lemma~\ref{lemma_sigma_2}}  we can split $\det(\I+Z_m)$ as follows: 
\begin{equation}\label{split}
\det(\I + Z_m)
=1 + \sigma_1(Z_m) + \sigma_2(Z_m) + \sigma_3(Z_m),
\end{equation}
where we repeat
\[
\sigma_1(Z_m) = \trace(Z_m), \quad
\sigma_2(Z_m) = \frac{1}{2} \left( \trace^2(Z_m)-\trace(Z_m^2)\right),\quad
\sigma_3(Z_m) = \det(Z_m).
\]

Thus
\begin{equation}\label{torsion_2m_twist}
{\mathbb T}^{J(2, 2m)}_{\lambda}(\rho)
=
\frac{\tau_0}{s+s^{-1}-2}
\left.{
    \frac{1}{2} 
    \frac{d^2}{dt^2}
        \left( 
            1 + \sigma_1(Z_m) + \sigma_2(Z_m) + \sigma_3(Z_m)
        \right)
    }\right|_{t=1}.
\end{equation}

With the \lq\lq splitting\rq\rq\, of  ${\mathbb T}^{J(2, 2m)}_{\lambda}(\rho)$ given in Eq.~(\ref{torsion_2m_twist}) in mind, we compute separately each \lq\lq second differential" of the $\sigma_i(Z_m)\,(i=1,2,3)$ to obtain the non--abelian Reidemeister torsion of $J(2, 2m)$.
\medskip

\noindent {\it The \lq\lq second differential" of $\sigma_3(Z_m)$.}
We concentrate first on  the $\sigma_3(Z_m)$-part of Eq.~(\ref{torsion_2m_twist}), which is the easier term to compute and correspond to the \lq\lq second differential" of $\sigma_3(Z_m)$.
\begin{claim}\label{claim:second_diff_sigma_3}
We have: 
\[
\frac{1}{2} \left.\frac{d^2}{dt^2} (\sigma_3(Z_m))\right|_{t=1} = \lim_{t \to 1}\frac{\sigma_3(Z_m)}{(t-1)^2}
=(2-s-s^{-1})^2 \cdot m t_m^2.
\]
\end{claim}
\begin{proof}[Proof of Claim \ref{claim:second_diff_sigma_3}]
By definition $\sigma_3(Z_m)=\det(Z_m)$, thus
\begin{align}
\sigma_3(Z_m)
    &=
    \det(
        (\I-tX)S_m(W^{-1})(t^{-1}X^{-1} - X^{-1}Y)
            ) \nonumber \\
    &=
    \det(\I-tX)
    \det(S_m(W^{-1}))
    \det(t^{-1}X^{-1}Y)
    \det(Y^{-1}-t\I) \nonumber \\
    &=
    t^{-3}(t-1)^2(1-ts)(1-ts^{-1})(t-s)(t-s^{-1})
    \det(S_m(W^{-1})). \label{sigma_3}
\end{align}
{Dividing Eq.~(\ref{sigma_3}) by $(t-1)^2$ and taking the limit} when $t$ goes to $1$, we thus obtain: 
\begin{equation*}\label{sigma_31}
\lim_{t \to 1}
\frac{\sigma_3(Z_m)}{(t-1)^2}
=
(2-s-s^{-1})^2\det(S_m(W^{-1})).
\end{equation*}
Note, with Eq.~(\ref{ev}) in mind, that $\xi_{\pm}^2$ and $1$ are the eigenvalues of $W^{-1} = Ad_{\rho(w)^{-1}}$. It thus follows that 
% \begin{align}
% \det(S_m(W^{-1})) \nonumber
% &=
% m
% \frac{(1-\xi_{+}^{2m})(1-\xi_{-}^{2m})}{(1-\xi_{+}^2)(1-\xi_{-}^2)} \nonumber\\
% &=
% m
% \frac{(\xi_{-}^m - \xi_{+}^m)(\xi_{+}^m - \xi_{-}^m)}{(\xi_{-} - \xi_{+})(\xi_{+} - \xi_{-})} \nonumber\\
% &=
% mt_m^2. \label{det_sum_W}
% \end{align}
\begin{align*}
\det(S_m(W^{-1}))
&=
m
\frac{(1-\xi_{+}^{2m})(1-\xi_{-}^{2m})}{(1-\xi_{+}^2)(1-\xi_{-}^2)} \\
&=
m
\frac{(\xi_{-}^m - \xi_{+}^m)(\xi_{+}^m - \xi_{-}^m)}{(\xi_{-} - \xi_{+})(\xi_{+} - \xi_{-})} \\
&=
mt_m^2. 
\end{align*}
\end{proof}

%If we substitute Equations~(\ref{sigma_3}) \&~(\ref{det_sum_W}) into Equation~(\ref{torsion_2m_twist}), 
If we substitute the result of Claim \ref{claim:second_diff_sigma_3} into Eq.~(\ref{torsion_2m_twist}), 
we obtain
%
%\begin{equation*}
%{\mathbb T}^{J(2, 2m)}_{\lambda}(\rho) = \frac{\tau_0}{s+s^{-1}-2}
%\left[
%   \lim_{t \to 1}
%   \frac{1+ \sigma_1(Z_m) + \sigma_2(Z_m)}{(t-1)^2}
%
%   + (s+s^{-1} - 2)^2 \cdot mt_m^2
%\right].
%\end{equation*}
%These expression can be easily written again as
\begin{equation}\label{sec_dif}
{\mathbb T}^{J(2, 2m)}_{\lambda}(\rho) = \frac{\tau_0}{s+s^{-1}-2}
\left[
    \left.
    \frac{1}{2} 
    \frac{d^2}{dt^2}
        \left( 
            \sigma_1(Z_m) + \sigma_2(Z_m)
        \right)
    \right|_{t=1}
+ (s+s^{-1} - 2)^2 \cdot mt_m^2
\right]. 
\end{equation}
%
%We will find non--vanishing terms in terms of second differential in $\sigma_1(Z_m)$ and $\sigma_2(Z_m)$.

\medskip

\noindent{\it The \lq\lq second differentials" of $\sigma_1(Z_m)$ and $\sigma_2(Z_m)$.}
We now focus on 
the \lq\lq second differentials" of  $\sigma_1(Z_m)$ and $\sigma_2(Z_m)$.
If we let 
$$\tilde Z_m = (t^{-1}X^{-1}-X^{-1}Y)(\I-tX)S_m(W^{-1}),$$
 then it follows from the definitions of $\sigma_1$ and $\sigma_2$ that 
\[
\sigma_1(Z_m) = \sigma_1(\tilde Z_m), \quad \sigma_2(Z_m) = \sigma_2(\tilde Z_m).
\]
We use $\sigma_1(\tilde Z_m)$ and $\sigma_2(\tilde Z_m)$ instead of 
$\sigma_1(Z_m)$ and $\sigma_2(Z_m)$ for our calculations.

\begin{claim}\label{claim:sec_diffs}
We have:
% second differential of sigma_1
\begin{align}
\frac{1}{2} \left. \frac{d^2}{dt^2}\sigma_1(\tilde Z_m) \right|_{t=1}
&= 
\trace(X^{-1}S_m(W^{-1})),  \label{sec_dif_sigma_1}\\
% second differential of sigam_2
\left.\frac{1}{2}\frac{d^2}{dt^2} \sigma_2(\tilde Z_m) \right|_{t=1} 
&= 3\sigma_2(X^{-1}S_m(W^{-1})) + \sigma_2(YS_m(W^{-1})W^{-1}) \label{sec_dif_sigma_2}\\
&\quad -\trace\left(X^{-1}S_m(W^{-1})\right)\trace\left((\I+X^{-1}Y)S_m(W^{-1})\right) \nonumber\\
&\quad + \trace\left(X^{-1}S_m(W^{-1})(\I+X^{-1}Y)S_m(W^{-1})\right).\nonumber 
\end{align}
\end{claim}
\begin{proof}[Proof of Claim \ref{claim:sec_diffs}]
Since $\sigma_1(\tilde Z_m)$ is the trace of $\tilde Z_m$, the
only term which remains after taking the \lq\lq second differential" at $t = 1$ is $\left. \frac{d^2}{dt^2}\frac{1}{t} X^{-1} S_m(W^{-1}) \right|_{t=1}$.

Now we consider $\sigma_2(\tilde Z_m)$.
From the definition of $\sigma_2(\tilde Z_m)$,
we have 
\begin{align*}
\frac{1}{2} 
\left. \frac{d^2}{dt^2}
\sigma_2(\tilde Z_m)
\right|_{t=1} 
&=
    \frac{1}{2}
        \left. \frac{d^2}{dt^2}
        \frac{1}{2}
        \left( 
            \trace^2(\tilde Z_m) -\trace ({\tilde Z_m}^2) 
        \right)
        \right|_{t=1}\\ 
&=
    \frac{1}{4}
    \left[  
    \left. \frac{d^2}{dt^2}\trace^2(\tilde Z_m) \right|_{t=1}
     - \left.\frac{d^2}{dt^2} \trace({\tilde Z_m}^2)\right|_{t=1}
    \right].
\end{align*}
In $\trace({\tilde Z_m}^2)$, 
the following three terms
are the terms which remain after taking the second differential at $t=1$:
\begin{eqnarray*}
&\left. \frac{d^2}{dt^2}t^{-2}\trace\left( (X^{-1}S_m(W^{-1}))^2 \right)\right|_{t=1},\\
&\left. \frac{d^2}{dt^2}-2t^{-1}\trace\left( X^{-1}S_m(W^{-1}) (\I+X^{-1}Y)S_m(W^{-1}) \right)\right|_{t=1},\\
&\left. \frac{d^2}{dt^2}t^2 \trace\left( (YS_m(W^{-1})W^{-1})^2 \right)\right|_{t=1}.
\end{eqnarray*}
Hence
\begin{align*} 
\left.\frac{1}{2}\frac{d^2}{dt^2} \sigma_2(\tilde Z_m) \right|_{t=1}
&=\frac{3}{2}
    \left[
        \trace^2\left(X^{-1}S_m(W^{-1})\right) - \trace\left( (X^{-1}S_m(W^{-1}))^2 \right) 
    \right] \\
&\quad -\trace\left(X^{-1}S_m(W^{-1})\right)\trace\left((\I+X^{-1}Y)S_m(W^{-1})\right) \\
&\quad + \trace\left(X^{-1}S_m(W^{-1})(\I+X^{-1}Y)S_m(W^{-1})\right)\\
&\quad +\frac{1}{2}
    \left[
        \trace^2\left(YS_m(W^{-1})W^{-1}\right) - \trace\left((YS_m(W^{-1})W^{-1})^2\right)
    \right].
\end{align*}
\end{proof}

If we substitute Eqs.~(\ref{sec_dif_sigma_1}) \&~(\ref{sec_dif_sigma_2}) of Claim~\ref{claim:sec_diffs}
into Eq.~(\ref{sec_dif}), then we obtain the following formula for ${\mathbb T}^{J(2, 2m)}_{\lambda} (\rho)$.

\begin{claim}\label{torsionJ(2,2m)}
The non--abelian Reidemeister torsion for $J(2, 2m)$ satisfies the following formula:
\begin{align} \label{matrix_rep_torsion_J(2,2m)}
{\mathbb T}^{J(2, 2m)}_{\lambda} (\rho)
&=
\frac{\tau_0}{s+s^{-1}-2}
\left[
    \trace\left(X^{-1}S_m(W^{-1})\right)
    +3\sigma_2(X^{-1}S_m(W^{-1})) \right. \\
    &\quad + \sigma_2(YS_m(W^{-1})W^{-1}) \nonumber \\ 
    &\quad -\trace\left(X^{-1}S_m(W^{-1})\right)\trace\left( S_m(W^{-1})\right) \nonumber \\
    &\quad + \trace\left(X^{-1}S_m(W^{-1})S_m(W^{-1})\right) \nonumber \\
 &\quad- \trace\left(X^{-1}S_m(W^{-1})\right)\trace\left( X^{-1}YS_m(W^{-1}) \right) \nonumber \\
 &\quad + \trace\left(X^{-1}S_m(W^{-1}) X^{-1}YS_m(W^{-1})\right)  
 \left.+ (s+s^{-1} - 2)^2 \cdot mt_m^2 \right]. \nonumber
\end{align}
\end{claim}

%%%%%%%%%%%%%%%%%%%%%%%%%%%%%%%%%%
\noindent{\it More explicit descriptions.}
%%%%%%%%%%%%%%%%%%%%%%%%%%%%%%%%%
To find more explicit expression of ${\mathbb T}^{J(2, 2m)}_{\lambda}(\rho)$,
we change our basis of $\sll$ in order to diagonalize the matrix {$\rho(w)$}.

The $\SL$-matrix $\rho(w)$ can be diagonalized by 
\[
p=
\left(
\begin{array}{cc}
u+1-s^{-1} & u+1-s^{-1} \\
1-su-\xi_{+} & 1-su- \xi_{-}
\end{array}
\right).
\]
Explicitly, $p^{-1} \rho(w) p$ is the diagonal matrix $\mathrm{diag}(\xi_+, \xi_-)$.

Set $a = 1-su-\xi_{+}$ and $b = 1-su-\xi_{-}$.
With respect to the basis $\{E, H, F\}$ of $\sll$, the matrix of the adjoint action of $p$ becomes as follows: 
\[
P = Ad_p = \frac{1}{a-b}
\left(
\begin{array}{ccc}
-c                         & 2c                       & c \\
 a                         & -(a+b)                   & -b \\
 a^2 / c                   & -2ab / c                 & - b^2 / c
\end{array}
\right)
\]
where $c = u+1-s^{-1}$ is defined in Subsection~\ref{statement}.

%\begin{remark}
Note that 
the matrix $P^{-1} W P$ is the diagonal matrix 
$\mathrm{diag}(\xi_{+}^2, 1, \xi_{-}^2)$. Here we repeat that $W = Ad_{\rho(w)}$.
%\end{remark}

Set 
$$\tilde X = P^{-1}XP,\quad  \tilde Y = P^{-1}YP \text{ and } \tilde W = P^{-1}WP.$$ 
Since we have $P^{-1}S_m(W^{-1}) P = S_m( \tilde W^{-1})$, 
the matrix $P^{-1} S_m (W^{-1}) P$ is the following diagonal matrix 
$$P^{-1} S_m (W^{-1}) P = \mathrm{diag}(\xi_{-}^{m-1} t_m, m, \xi_{+}^{m-1} t_m).$$

Moreover as $\trace(X^{-1}S_m(W^{-1})) = \trace(\tilde X^{-1} S_m(\tilde W^{-1}))$, we have the following claim.
\begin{claim}\label{trace_X^{-1}S_m}
We have
\begin{align*}
\lefteqn{
\trace(X^{-1}S_m(W^{-1})) }
& \\
&=
\frac{1}{(a-b)^2}
\left(
\frac{\beta_1^2}{s} \xi_{-}^{m-1} t_m 
+ \left((a-b)^2(s+s^{-1}+1)-\frac{\beta_1^2}{s} - \frac{\alpha_1^2}{s}\right)m
+\frac{\alpha_1^2}{s} \xi_{+}^{m-1} t_m
\right).
\end{align*}
\end{claim}
\begin{proof}[Proof of Claim \ref{trace_X^{-1}S_m}]
By a direct computation, we obtain that the $(1, 1)$-component of the matrix $\tilde X^{-1}$ is equal to $\beta_1^2 / (s(a-b)^2)$ and its 
$(3, 3)$-component is equal to $\alpha_1^2 / (s(a-b)^2)$. 
We can also find the $(2,2)$-component of $\tilde X^{-1}$ from $\trace(\tilde X^{-1}) = s+s^{-1}+1$.
\end{proof}

Now, we compute $\sigma_2(XS_m(W^{-1}))$ and $\sigma_2(YS_m(W^{-1})W^{-1})$
from Lemma \ref{lemma_sigma_2} as follows.
\begin{claim}\label{sigma_2_result}
The following equalities hold: 
\begin{align*}
\lefteqn{
  (1) \quad
  \sigma_2(X^{-1}S_m(W^{-1}))
}
& \\
&\quad=
  \frac{1}{(a-b)^2}
  \left(
    \frac{\beta^2_1}{s} \xi_{-}^{m-1}m t_m 
    + \left((a-b)^2(s+s^{-1}+1) - \frac{\alpha^2_1}{s} - \frac{\beta^2_1}{s} \right) t_m^2 + \frac{\alpha^2_1}{s} \xi_{+}^{m-1} m t_m
  \right).\\
\lefteqn{
  (2) \quad
  \sigma_2(YS_m(W^{-1})W^{-1})
}
& \\
&\quad=
  \frac{1}{(a-b)^2}
  \left(
    \frac{\alpha^2_1}{s} \xi_{-}^{m+1}m t_m 
    + \left((a-b)^2(s+s^{-1}+1) - \frac{\alpha^2_1}{s} - \frac{\beta^2_1}{s} \right) t_m^2
    + \frac{\beta^2_1}{s} \xi_{+}^{m+1} m t_m
  \right).
\end{align*}
\end{claim}
\begin{proof}[Proof of Claim \ref{sigma_2_result}]
Using Lemma \ref{lemma_sigma_2} and because $\det(X)=1, \det(Y) = 1$ we have:
\[
\sigma_2(X^{-1}S_m(W^{-1})) = \trace(X S_m(W^{-1})^{-1})\cdot \det(S_m(W^{-1}))
\] 
and 
\[
\sigma_2(YS_m(W^{-1}) W^{-1}) = \trace(Y^{-1} W S_m(W^{-1})^{-1})\cdot \det(S_m(W^{-1})).
\]
We obtain the above formulas by computing the traces using $\tilde X, \tilde Y$ and $S_m(\tilde W^{-1})$ as in Claim~\ref{trace_X^{-1}S_m}.
\end{proof}

Finally we calculate the other two terms 
which are of the following form: $-\trace(A)\trace(B) + \trace(AB)$.
\begin{claim}\label{trAB-trAtrB}
We have:
\begin{align*}
(1)\quad\lefteqn{
-\trace(X^{-1}S_m(W^{-1}))\trace(S_m(W^{-1}))
+\trace(X^{-1}S_m(W^{-1})S_m(W^{-1}))
}
& \\
&=
- \frac{1}{(a-b)^2}
\left[
\left(
(a-b)^2(s+\frac{1}{s}+1) -\frac{\alpha^2_1}{s} 
\right)\xi_{-}^{m-1} m t_m
+\left(
\frac{\alpha^2_1}{s} + \frac{\beta^2_1}{s} 
\right)t^2_m \right.\\
& \left.\quad+\left(
(a-b)^2(s+\frac{1}{s}+1) -\frac{\beta^2_1}{s} 
\right)\xi_{+}^{m-1} m t_m.
\right]
\end{align*}
(2) If we set $\tilde X^{-1}={(a_{i, j})}_{1\leq i, j \leq 3},$ and  $\tilde X^{-1} \tilde Y={(b_{i, j})}_{1 \leq i, j \leq 3}$,
then we have
\begin{align*}
&
-\trace(X^{-1}S_m(W^{-1}))\trace(X^{-1} Y S_m(W^{-1}))
+\trace(X^{-1}S_m(W^{-1})X^{-1} Y S_m(W^{-1}))
\\
&=
\left(
    \left|
    \begin{array}{cc}
    a_{1, 1} & a_{1, 3} \\
    b_{3, 1} & b_{3, 3}
    \end{array}
    \right|
    +
    \left|
    \begin{array}{cc}
    b_{1, 1} & b_{1, 3} \\
    a_{3, 1} & a_{3, 3} 
    \end{array}
    \right|
\right)
t^2_m \\
& \quad -
\left(
    \left|
    \begin{array}{cc}
    a_{2, 2} & a_{2, 3} \\
    b_{3, 2} & b_{3, 3}
    \end{array}
    \right|
    +
    \left|
    \begin{array}{cc}
    b_{2, 2} & b_{2, 3} \\
    a_{3, 2} & a_{3, 3} 
    \end{array}
    \right|
\right)m 
\xi_{+}^{m-1} t_m \\
& \quad -
\left(
    \left|
    \begin{array}{cc}
    a_{1, 1} & a_{1, 2} \\
    b_{2, 1} & b_{2, 2}
    \end{array}
    \right|
    +
    \left|
    \begin{array}{cc}
    b_{1, 1} & b_{1, 2} \\
    a_{2, 1} & a_{2, 2} 
    \end{array}
    \right|
\right)
m \xi_{-}^{m-1}t_m.
\end{align*}
\end{claim}
\begin{proof}[Proof of Claim \ref{trAB-trAtrB}]
Item (1) follows from above results and Item (2) follows from Lem\-ma~\ref{lemma_sigma_2}.
\end{proof}

\begin{remark}
 The  matrices $\tilde X^{-1}$ and $\tilde X^{-1} \tilde Y$ are described explicitly as follows.
\[
\tilde X ^{-1}
=
\frac{1}{(a-b)^2}
\left(
\begin{array}{ccc}
\frac{\beta^2_1}{s}  & -\frac{2\beta_1 \beta_2}{s} & -\frac{\beta^2_2}{s} \\
\frac{\alpha_2 \beta_1}{s} & (a-b)^2(s+\frac{1}{s}+1)-\frac{\alpha^2_1+\beta^2_1}{s} & -\frac{\alpha_1 \beta_2}{s} \\
-\frac{\alpha_2^2}{s} & \frac{2\alpha_1 \alpha_2}{s} & \frac{\alpha^2_1}{s}
\end{array}
\right),
\]
\begin{align*}
\lefteqn{
(a-b)^2 \tilde X^{-1} \tilde Y
} & \\
&=
\left(
\begin{array}{lll}
(c(1-\xi_{-})+\frac{\beta_1}{s})^2 & -2(c(1-\xi_{-})+\frac{\beta_1}{s})(c(1-\xi_{-})+\frac{\beta_2}{s}) & -(c(1-\xi_{-})+\frac{\beta_2}{s})^2 \\
\begin{array}{l}
(c(1-\xi_{+})+\frac{\alpha_2}{s})\\
\quad\cdot(c(1-\xi_{-})+\frac{\beta_1}{s})
\end{array}
&
\begin{array}{l}
 (a-b)^2(u^2+4u+3) \\
 \quad - (c(1-\xi_{+})+\frac{\alpha_1}{s})^2 \\
 \quad -(c(1-\xi_{-})+\frac{\beta_1}{s})^2
\end{array}
 &
 \begin{array}{l}
-(c(1-\xi_{+})+\frac{\alpha_1}{s})\\
 \quad \cdot (c(1-\xi_{-})+\frac{\beta_2}{s})
\end{array}
\\
-(c(1-\xi_{+})+\frac{\alpha_2}{s})^2 & 2(c(1-\xi_{+})+\frac{\alpha_1}{s})(c(1-\xi_{+})+\frac{\alpha_2}{s}) &(c(1-\xi_{+})+\frac{\alpha_1}{s})^2 
\end{array}
\right).
\end{align*}
\end{remark}

One can observe that $(a-b)^2$ is equal to $(\xi_{+} - \xi_{-})^2$. Thus, 
if we substitute the results given in
Claims \ref{trace_X^{-1}S_m}, \ref{sigma_2_result} \& \ref{trAB-trAtrB}
into Eq.~(\ref{matrix_rep_torsion_J(2,2m)}), then we obtain Eq.~(\ref{m_pos}) of Theorem~\ref{theorem_J(2,2m)}. This achieves the first part of the proof of Theorem~\ref{theorem_J(2,2m)}. 

The computation of ${\mathbb T}^{J(2, -2m)}_{\lambda} (\rho)$, for $m>0$, is completely similar and has the following expression:
\begin{align}
{\mathbb T}^{J(2, -2m)}_{\lambda}(\rho)
&=
-\tau_0
\lim_{t \to 1}
\frac{\det ( \Phi(\frac{\partial}{\partial x}w^{-m}  x w^m y^{-1})) }
            {(t-1)\det ( \Phi(y-1))} \nonumber \\
&=
\frac{\tau_0}{s+s^{-1} - 2}
\lim_{t \to 1}
\frac{\det(\I+Z_{-m})}{(t-1)^2}. \label{eqn:torsion_J(2,-2m)}
\end{align}
Here the matrix $Z_{-m}$ is given by 
$(\I-tX)S_m(W)(YX^{-1}- t^{-1}YX^{-1}Y^{-1})$.
The right hand side of Eq.~(\ref{eqn:torsion_J(2,-2m)}) is given by
\begin{align}\label{eqn:torsion_J(2,-2m)_tr}
&\frac{\tau_0}{s+s^{-1}-2}
\left[
  -\trace(X^{-1}S_m(W)W)
   +3\sigma_2(WX^{-1}S_m(W)) + \sigma_2(YS_m(W))
\right. \\
&\quad 
-\trace(X^{-1}S_m(W)W)\trace((\I+X^{-1}Y)S_m(W)W)  \nonumber\\
&\quad 
+ \trace(X^{-1}S_m(W)W(\I+X^{-1}Y)S_m(W)W)  \nonumber \\
&\quad \left.- (s+s^{-1} - 2)^2 \cdot mt_m^2 \right]. \nonumber
\end{align}
Each term in Eq.~(\ref{eqn:torsion_J(2,-2m)_tr}) can be computed similarly as in Claims~\ref{trace_X^{-1}S_m},~\ref{sigma_2_result} \&~\ref{trAB-trAtrB} which give the second item of Theorem~\ref{theorem_J(2,2m)}.

%%%%%%%%%%%%%%%%%%%%%%%%%%%%%%%%%%%%%%%%%%%%%%%%
\subsection{Examples}
\label{Sexamples}
%%%%%%%%%%%%%%%%%%%%%%%%%%%%%%%%%%%%%%%%%%%%%%%%
As an illustration of our main Theorem~\ref{theorem_J(2,2m)}, we explicitly compute ${\mathbb T}^{K}_{\lambda}$ on the character variety $X^{\mathrm{irr}}_\lambda(\Pi(K))$ for the following four examples: the trefoil knot $3_1 = J(2,2)$, $5_2 = J(2,4)$, the figure eight knot $4_1 = J(2,-2)$, and $6_1 = J(2,-4)$ respectively. The computer program for the computation is written on the free computer algebra system Maxima~\cite{Max}.

\begin{itemize}
\item[(1)]
    For the trefoil knot $3_1 = J(2, 2)$, the Riley polynomial is given by 
    $$\phi_{J(2, 2)}(s, u) = -1 + s+s^{-1} - u.$$
    The computation of the non--abelian Reidemeister torsion for $J(2, 2)$ is expressed as follows.
    \begin{align*}
        {\mathbb T}^{J(2, 2)}_{\lambda}(\rho_{\sqrt{s}, u}) 
            &= \frac{\tau_0}{s + s^{-1} -2} \left(-3(s+s^{-1}-2)\right) \\
            &= -3 \tau_0.
    \end{align*}
    This coincides with the inverse of the result~\cite[Subsection~6.1]{JDFibre} (see Remark~\ref{rmkinv}).
\item[(2)]
    For $5_2 = J(2, 4)$, the Riley polynomial is given by 
    \begin{align*}
    \phi_{J(2, 4)}(s, u) =-3+2(s+s^{-1}) &+\left(-4+3(s+s^{-1})-(s+s^{-1})^2\right)u \\ &+ \left(-3+2(s+s^{-1})\right) u^2 - u^3.
    \end{align*}
    The non--abelian Reidemeister torsion for $5_2 = J(2, 4)$ is expressed as follows.
    \begin{align*}
        {\mathbb T}^{J(2, 4)}_{\lambda}(\rho_{\sqrt{s}, u}) 
            &= \frac{\tau_0}{s + s^{-1} -2}  
                   [ -2 + 21(s+s^{-1})-10(s+s^{-1})^2 \\
            &\quad
                    +\{-2+15(s+s^{-1})-17(s+s^{-1})^2+5(s+s^{-1})^3\}u \\
            &\quad
                    +\{6+7(s+s^{-1})-5(s+s^{-1})^2\}u^2 ]\\
            &=\tau_0
              \left(-10(s+s^{-1})+1
                +\left(5(s+s^{-1})^2-7(s+s^{-1})+1\right)u \right.\\
            &\qquad\qquad
               \left.+\left(-5(s+s^{-1})-3\right)u^2
               \right).
    \end{align*}
\item[(3)]
    For the figure eight knot $4_1 = J(2, -2)$, the Riley polynomial is given by  $$\phi_{J(2, -2)}(s, u)=(3-s-s^{-1})(u+1)+u^2.
    $$
    The computation of the non--abelian Reidemeister torsion for $J(2, -2)$ is expressed as follows.
    \begin{align*}
        {\mathbb T}^{J(2, -2)}_{\lambda}(\rho_{\sqrt{s}, u}) 
            &= \frac{\tau_0}{s+s^{-1}-2} \left(-2(s+s^{-1})+1\right)\left(s+s^{-1} -2\right) \\
            &= \tau_0 \left(-2(s+s^{-1})+1\right).
    \end{align*}
    This coincides with the inverse of the result~\cite[Subsection 6.3]{JDFibre} (see Remark~\ref{rmkinv}) in which the torsion is expressed as $\pm \sqrt{17 + 4I_{\lambda}} $. 
    
    Since the longitude $\lambda$ is equal to $[y,x^{-1}][x, y^{-1}]$, one has  
    $$I_{\lambda} = -2 - (s+s^{-1}) + s^2+s^{-2}.$$
 Thus, up to sign, we have 
$$\sqrt{17 + 4I_{\lambda}}= 2(s+s^{-1})-1.$$
\item[(4)]
    For $6_1 = J(2, -4)$, the Riley polynomial is given by 
\[
\begin{split}
\phi_{J(2, -4)}(s, u)= 5 &- 2(s+s^{-1})+\left(12+(s+s^{-1})^2 -7(s+s^{-1})\right)u + \\ & \left(11+(s+s^{-1})^2-6(s+s^{-1})\right)u^2+\left(5-2(s+s^{-1})\right)u^3+u^4.
\end{split}
\]
    The non--abelian Reidemeister torsion for $J(2, -4)$ is expressed as follows.
    \begin{align*}
    \lefteqn{
        {\mathbb T}^{J(2, -4)}_{\lambda}(\rho_{\sqrt{s}, u})} &\\
            &= \frac{\tau_0}{s + s^{-1} -2}  
                   [-14 - 13(s+s^{-1}) +26(s+s^{-1})^2 -8(s+s^{-1})^3\\
            &\quad
                    +\{-8 - 34(s+s^{-1}) +49(s+s^{-1})^2 - 23(s+s^{-1})^3 + 4(s+s^{-1})^4\}u \\
            &\quad
                    +\{-2 - 31(s+s^{-1}) + 32(s+s^{-1})^2 - 8(s+s^{-1})^3\}u^2 \\
            &\quad 
                    +\{2 -9(s+s^{-1}) +4(s+s^{-1})^2\}u^3 ]\\
            &=\tau_0\left[ \left(-8(s+s^{-1})^2+10(s+s^{-1})+7\right) \right.\\
            &\qquad +\left(4(s+s^{-1})^3-15(s+s^{-1})^2+19(s+s^{-1})+4\right)u\\
            &\qquad \left.+\left(-8(s+s^{-1})^2+16(s+s^{-1})+1\right)u^2
                   +\left(4(s+s^{-1})-1\right)u^3 \right].
    \end{align*}
\end{itemize}
%%%%%%%%%%%%%%%%%%%%%%%%%%%%%%%%%%%%%%%%%%%%%%%%
\subsection{Twisted Reidemeister torsion at the holonomy representation.}
\label{subsection:results_holonomy}
%%%%%%%%%%%%%%%%%%%%%%%%%%%%%%%%%%%%%%%%%%%%%%%%
In this section, we consider the non--abelian Reidemeister torsion for hyperbolic twist knots at holonomy representations.
Formulas of the non--abelian Reidemeister torsion associated to twist knots are complicated.
But we see here that formulas for the non--abelian Reidemeister torsion at holonomy representations are simpler.

Every twist knots except the trefoil knot are hyperbolic. 
It is well known that 
an exterior of a hyperbolic knot admits at most a complete hyperbolic structure and 
this hyperbolic structure determines the holonomy representation of the knot group (see Section~\ref{sechyp}). With Fact~\ref{holonomy} in mind, we know that such lifts are $\lambda$-regular representations.

\begin{remark}
If we substitute $s=1$ into Riley's polynomial $\phi_{J(2, 2m)} (s, u)$ given in Eq.~(\ref{eqn:Riley_poly}), then
\[
\phi_{J(2, 2m)} (1, u)= (1-u) t_m - t_{m-1} = 
(1-u) \frac{\xi_{+}^m - \xi_{-}^m}{\xi_{+} - \xi_{-}} - \frac{\xi_{+}^{m-1}-\xi_{-}^{m-1}}{\xi_{+} - \xi_{-}}.
\] 
\end{remark}
%
%\begin{remark}
%It is easy to express the root $u_0$ of Riley's equation $\phi(1,u)$ corresponding to the holonomy in terms of the cusp shape $\cusp$ of the knot:
%\begin{equation}\label{cuspshapeTK}
%u_0 = \frac{4}{\cusp-2}.
%\end{equation}
%\end{remark}

%
The $\SL$-representations which lifts the holonomy representation correspond 
to roots of Riley's equation $\phi_{J(2, 2m)}(1,u)=0$. We let $\rho_u$ denote such representations.

We are now ready to give a formula for the non--abelian Reidemeister torsion of $J(2,2m)$ at the holonomy representation in terms of $u$, which is an algebraic number, defined to be the root of  a computable algebraic equation. 

\begin{theorem}\label{R-torsionTK}
Let $m>0$ and $u$ denote one of the two complex conjugate roots of Riley's equation $\phi_{J(2, 2m)}(1,u) = 0$ (resp. $\phi_{J(2, -2m)}(1,u) = 0$) corresponding to holonomy representations, 
then 
\begin{enumerate}
\item the non--abelian Reidemeister torsion ${\mathbb T}^{J(2, 2m)}_{\lambda}(\rho_u)$ satisfies the following
{formula}:
\begin{align*}
{\mathbb T}^{J(2, 2m)}_{\lambda}(\rho_u)
&= 
  \frac{-\tau_0}{u^2+4}
  \left[
    \left(4+m(u^2-4u+8)\right) t_m (\xi_{+}^m + \xi_{-}^m) 
    \right. \nonumber \\
& \quad + 
    \left(t_m(\xi_{+}^{m-1}+\xi_{-}^{m-1})-1\right)(u^2-4)m \nonumber\\
& \left.\quad+
    (-5u^2-8u+4)t_m^2
  \right],
\end{align*}
\item similarly
the Reidemeister torsion ${\mathbb T}^{J(2, -2m)}_{\lambda}(\rho_u)$ satisfies the following
formula:
\begin{align*}
{\mathbb T}^{J(2, -2m)}_{\lambda}(\rho_u)
&= 
  \frac{-\tau_0}{u^2+4}
  \left[
    \left(-4+m(u^2-4u+8)\right) t_m (\xi_{+}^m + \xi_{-}^m) 
    \right. \nonumber \\
& \quad + 
    \left(t_m(\xi_{+}^{m+1}+\xi_{-}^{m+1}) +1\right)(u^2-4)m \nonumber\\
& \left.\quad+
    (-5u^2-8u+4)t_m^2
 \right].
\end{align*}
\end{enumerate}
\end{theorem}
\begin{remark}
Combining results of Theorem~\ref{R-torsionTK} with Eq.~(\ref{cuspshape}), the non--abelian Reidemeister torsion of a (hyperbolic) twist knot at the holonomy is expressed in terms of the cusp shape of the knot.
\end{remark}

\begin{proof}
First we make the computations in the case of $J(2, 2m)$, where $m > 0$. If we substitute $s = 1$ in Eq.~(\ref{Tlim}), then we obtain:
\begin{align}
{\mathbb T}^{J(2, 2m)}_{\lambda} (\rho)
&=
-\tau_0
\lim_{t \to 1}
\frac{\det(\I + Z_m)}{(t-1)^4} \nonumber \\
&=
\frac{-\tau_0}{24} \left.\frac{d^4}{dt^4}\det(\I+Z_m)\right|_{t=1}. \nonumber
\end{align}
Next, using the splitting of $\det(\I+Z_m)$ given in Eq.~(\ref{split}), we get:
\begin{equation}\label{diff_hol}
{\mathbb T}^{J(2, 2m)}_{\lambda} (\rho)
 = \frac{-\tau_0}{24} \left.\frac{d^4}{dt^4}\left[ 1^{} + \sigma_1(Z_m) + \sigma_2(Z_m) + \sigma_3(Z_m) \right]\right|_{t=1}. 
 \end{equation}
It follows from Eq.~(\ref{sigma_3}) that $(t-1)^6$ divides $\sigma_3(Z_m)$ in the case of $s = 1$.
Hence the term $\left.\frac{d^4}{dt^4} \sigma_3(Z_m)\right|_{t=1}$ in Eq.~(\ref{diff_hol}) vanishes.
By degrees of $t$ in $\sigma_1(Z_m)$ and $\sigma_2(Z_m)$, 
we obtain the following equation from direct computations of the above differentials:
\begin{align}
{\mathbb T}^{J(2, 2m)}_{\lambda} (\rho)
&=
-\tau_0[
\trace(X^{-1}S_m(W^{-1}))+5\sigma_2(X^{-1}S_m(W^{-1})) \label{tor_hol_matrix_pos}\\
& \quad
-\trace(X^{-1}S_m(W^{-1})) \trace((\I+X^{-1}Y)S_m(W^{-1})) \nonumber \\
&\quad
+\trace(X^{-1}S_m(W^{-1})(\I+X^{-1}Y)S_m(W^{-1})) ] \nonumber.
\end{align}
Note that $\xi_{+}  + \xi_{-} = u^2+2$ and $(\xi_{+} - \xi_{-})^2 = u^2(u^2+4)$.
By Claims \ref{trace_X^{-1}S_m}, \ref{sigma_2_result} \& \ref{trAB-trAtrB} 
we have
\begin{align*}
\trace(X^{-1}S_m(W^{-1})) 
&=
\frac{1}{u^2+4}
\left[
4(\xi_{+}^m + \xi_{-}^m)t_m - (u^2-4)m
\right] \\
\sigma_2(X^{-1}S_m(W^{-1})) 
&=
\frac{1}{u^2+4}
\left[
4m(\xi_{+}^m + \xi_{-}^m)t_m - (u^2 - 4)t^2_m
\right]
\end{align*}
and 
\begin{align*}
\lefteqn{
-\trace(X^{-1}S_m(W^{-1}))\trace((\I+X^{-1}Y)S_m(W^{-1}))
+ \trace(X^{-1}S_m(W^{-1})(\I+X^{-1}Y)S_m(W^{-1}))
}
& \\
&=
-\frac{1}{u^2+4}
\left[
8(u+2)t^2_m 
-m(u+2)(u-6)(\xi_{+}^m + \xi_{-}^m)t_m
-m(u^2-4)(\xi_{+}^{m-1}+\xi_{-}^{m-1})t_m
\right].
\end{align*}
If we substitute these results into Eq.~(\ref{tor_hol_matrix_pos}), then
we obtain the wanted formula for $J(2, 2m)$.

Similarly, in the case of twist knots $J(2, -2m)$, $m>0$, from computations of differentials we have 
\begin{align}
{\mathbb T}^{J(2, -2m)}_{\lambda}(\rho_u)
&=
-\tau_0[
-\trace(X^{-1}S_m(W)W)+5\sigma_2(X^{-1}S_m(W)W) \label{tor_hol_matrix_neg}\\
& \quad
-\trace(X^{-1}S_m(W)W) \trace((\I+X^{-1}Y)S_m(W)W) \nonumber \\
&\quad
+\trace(X^{-1}S_m(W)W(\I+X^{-1}Y)S_m(W)W) ] \nonumber.
\end{align}
It follows from Claims \ref{trace_X^{-1}S_m}, \ref{sigma_2_result} \& \ref{trAB-trAtrB} that
\begin{align*}
\trace(X^{-1}S_m(W)W)
&=
\frac{1}{u^2+4}
\left[
4(\xi_{+}^m + \xi_{-}^m)t_m - (u^2-4)m
\right] \\
\sigma_2(X^{-1}S_m(W)W) 
&=
\frac{1}{u^2+4}
\left[
4m(\xi_{+}^m + \xi_{-}^m)t_m - (u^2 - 4)t^2_m
\right]
\end{align*}
and 
\begin{align*}
\lefteqn{
-\trace(X^{-1}S_m(W)W)\trace((\I+X^{-1}Y)S_m(W)W)
+\trace(X^{-1}S_m(W)W(\I+X^{-1}YS_m(W)W))
}& \\
&=
-\frac{1}{u^2+4}
\left[
8(u+2)t^2_m 
-m(u+2)(u-6)(\xi_{+}^m + \xi_{-}^m)t_m
-m(u^2-4)(\xi_{+}^{m+1} + \xi_{-}^{m+1})t_m
\right].
\end{align*}
If we substitute these results into Eq. (\ref{tor_hol_matrix_neg}) , then
we obtain the wanted formula for $J(2, -2m)$.
\end{proof}

\subsection{A remark on the asymptotic behavior of the non--abelian Reidemeister torsion at holonomy}
\label{ExTb}
%\subsection{Explicit examples}
%We consider again our list of examples of Subsection~\ref{Sexamples}.

We close this paper with some remarks on the behavior of the cusp shape and of the non--abelian Reidemeister torsion at the holonomy for twist knots. 

\begin{figure}[!htbp]
\begin{center}
\scalebox{.4}{\includegraphics{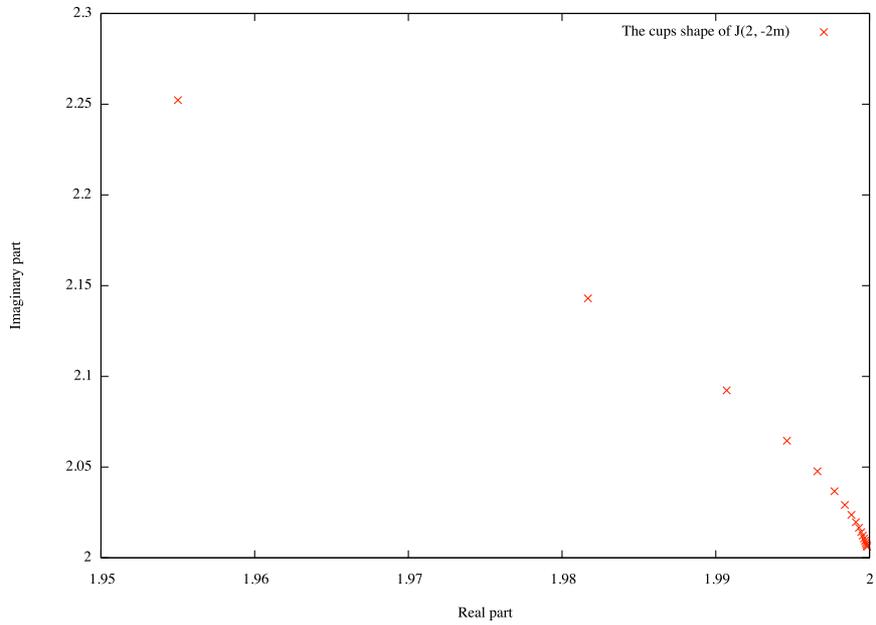}}
\end{center}
\caption{Graph of the cusp shape of $J(2,-2m)$ (from Snappea~\cite{snappea}).}
\label{graph:cusp_shapes}
\end{figure}

\begin{remark}[Behavior of the cusp shape]
In {the Notes}~\cite[p.~5.63]{TN}, Thurston explains that the sequence of {exteriors of the twist knots ${J(2, -2m)}$} converges to the exterior of the Whitehead link on Fig. \ref{fig:Whitehead} (link $5^2_1$ in Rolfsen's table~\cite{Rolfsen}). 
Note that, if the number of crossings $m$ increases to infinity, then the {sequence of cusp shapes} of the twist knots $J(2, \pm 2m)$ converges to $2+2i$, which is the common value of the cusp shapes of the Whitehead link, see the graph on Fig.~\ref{graph:cusp_shapes}. This result is a consequence of Dehn's hyperbolic surgery theorem.
\end{remark}

% Graph of Torsion and 0.0325*(2x+2)^3
%
\begin{figure}[!htbp]
\begin{center}
\scalebox{.4}{\includegraphics{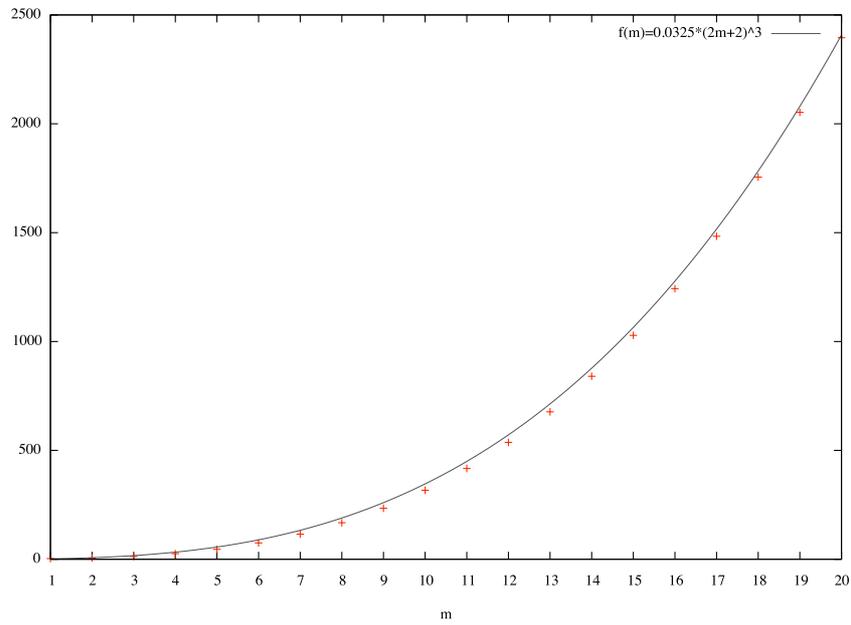}}
\end{center}
\caption{Graph of $|{\mathbb T}^{J(2, -2m)}_{\lambda}(\rho_0)|$ and $f(m)=C (\sharp J(2, -2m))^3$.}
\label{graph:compare}
\end{figure}

% % Graph of Torsion and 0.0325*(2x+2)^3
% %
% \begin{figure}[!htbp]
% \begin{center}
% \scalebox{.4}{\includegraphics{graphCompare}}
% \end{center}
% \caption{Graph of $|{\mathbb T}^{J(2, -2m)}_{\lambda}(\rho_0)|$ and $f(m)=C (\sharp J(2, -2m))^3$.}
% \label{graph:compare}
% \end{figure}

%\begin{remark}[Behavior of the torsion]
The graph on Fig.~\ref{graph:compare} gives the behavior of the sequence of the absolute value of  ${\mathbb T}^{J(2, -2m)}_{\lambda}(\rho_0)$ with respect to the number of crossings $\sharp J(2, -2m) = 2 + 2m$ of the knot. The order of growth can be deduced by a \lq\lq surgery argument\rq\rq\, using Item (5) of Remark~\ref{remarktwist} and the surgery formula for the Reidemeister torsion~\cite[Theorem 4.1]{Porti:1997}.
%\end{remark}

\begin{proposition}
The sequence $\left( {|{\mathbb T}^{J(2, -2m)}_{\lambda}(\rho_0)|} \right)_{m\geqslant 1}$ has the same behavior as the sequence $\left( {C (\sharp J(2, -2m))^3} \right)_{m\geqslant 1}$, for some constant $C$. 
\end{proposition}
\begin{proof}[Ideas of the Proof]
Item (5) of Remark~\ref{remarktwist} gives us that $E_{J(2, -2m)} = \mathcal{W}(1/m)$ is obtained by a surgery of slope $1/m$ on the trivial component of the Whitehead link $\mathcal{W}$. Let $V$ denote the glued solid torus and $\gamma$ its core. Using~\cite[Theorem 4.1 (iii) and Proposition 2.25]{Porti:1997} we have, up to sign:
\begin{equation}\label{surgeryF}
% {\mathbb T}^{J(2, -2m)}_{\lambda}(\rho_0) = \frac{{\mathbb T}^{W}_{(\lambda, \mu' {\lambda'}^{-m})}(\rho_0)}{\mathrm{TOR}(V; \sll_{\rho_0}, \gamma)}
{\mathbb T}^{J(2, -2m)}_{\lambda}(\rho_0) = {\mathbb T}^{\mathcal{W}}_{(\lambda, \mu' {\lambda'}^{-m})}(\rho_0) \cdot \mathrm{TOR}(V; \sll_{\rho_0}, \gamma)
\end{equation}
where ${\mathbb T}^{\mathcal{W}}_{(\lambda, \mu' {\lambda'}^{-m})}(\rho_0)$ stands for the $(Ad \circ \rho_0)$-non--abelian Reidemeister torsion of the Whitehead link exterior computed with respect to the bases of the twisted homology groups determined by the two curves $\lambda$, and $\mu' {\lambda'}^{-m}$ (see~\cite[Theorem 4.1]{Porti:1997}), and $\mathrm{TOR}(V; \sll_{\rho_0}, \gamma)$ stands for the $(Ad \circ \rho_0)$-non--abelian Reidemeister torsion of the solid torus $V$ computed with respect to its core $\gamma$. Here $\rho_0$ denotes the holonomy representation of the Whitehead link exterior.

To obtain the behavior of ${\mathbb T}^{J(2, -2m)}_{\lambda}(\rho_0)$ we estimate the two terms on the right--hand side of Eq.~(\ref{surgeryF}):
\begin{enumerate}
  \item Using~\cite[Proof of Theorem 4.17]{Porti:1997}, it is easy to see that $\mathrm{TOR}(V; \sll_{\rho_0}, \gamma)$ goes as $\frac{1}{4 \pi^2} m^2$ when $m$ goes to infinity.
  \item Using~\cite[Theorem 4.1 (ii)]{Porti:1997}, one can prove that 
  \[
  {\mathbb T}^{\mathcal{W}}_{(\lambda, \mu' {\lambda'}^{-m})}(\rho_0) = {\mathbb T}^{\mathcal{W}}_{(\lambda, \lambda')}(\rho_0) \cdot \left( \frac{1}{\cusp(\lambda, \mu)} - m\right),
  \] 
  where $\cusp = \cusp(\lambda, \mu)$ denotes the cusp shape of $J(2, -2m)$ (computed with respect to the usual meridian--longitude system). Thus, ${\mathbb T}^{\mathcal{W}}_{(\lambda, \mu' {\lambda'}^{-m})}(\rho_0)$ goes as $ {\mathbb T}^{\mathcal{W}}_{(\lambda, \lambda')}(\rho_0) \cdot m$ when $m$ goes to infinity. One can also prove that, at the holonomy, we have 
  \[
  {\mathbb T}^{\mathcal{W}}_{(\lambda, \lambda')}(\rho_0)  = 8(1+i).
  \]
\end{enumerate}
As a result, ${\mathbb T}^{J(2, -2m)}_{\lambda}(\rho_0)$ goes as $C\cdot m^3$ for some constant $C$.
\end{proof}

\clearpage
\appendix
\section{Program list and Tables}

\subsection{Program list for Maxima}
We give a program list 
in order to compute the twisted Reidemeister torsion 
for a given twist knot.
This program works on the free computer algebra system \emph{Maxima}~\cite{Max}.
The function $R(m)$ in the list computes the Riley polynomial of $J(2, 2m)$.
The function $T(m)$ computes  the twisted Reidemeister torsion for $J(2, 2m)$.
It gives a polynomial of $s$ and $u$ such that the top degree of $u$ is lower
than that in the Riley polynomial $\phi_{J(2, 2m)}(s, u)$.
Here we use Expressions~(\ref{matrix_rep_torsion_J(2,2m)}) \&~(\ref{eqn:torsion_J(2,-2m)_tr}) and the following remark
for computing the twisted Reidemeister torsion.
\begin{remark}
It follows from Equation~$(\ref{eqn:Riley_poly})$ that
the highest degree term of $u$ in $\phi_{J(2,2m)}(s, u)$ is equal to 
 $-u^{2m-1}$ (resp. $u^{2|m|}$) if $m>0$ (resp. $m<0$).
% \[
% \left\{
% \begin{array}{cc}
% -u^{2m-1} & m > 0,\\
% u^{2|m|} & m < 0.   
% \end{array}
% \right.
% \]
\end{remark}
\noindent{\bf Program list}
\begin{small}
\begin{verbatim}
load("nchrpl");/*We need this package for using mattrace*/
R(m):=block(/*function for calculating the Riley polynomial of J(2,2m)*/ 
        [/*w is the matrix of w=[y,x^{-1}]*/
         w:matrix([1-s*u,1/s-u-1],[-u+s*u*(u+1),(-u)/s+(u+1)^2]), 
         p],
        w:w^^m,
        p:w[1,1]+(1-s)*w[1,2],
        p:expand(p),
        return(p)
      );
T(m):=if integerp(m) then 
        if m=0 then "J(2,0) is unknot." else 
          block(
            [/*matrix for adjoint action of x*/
             X:matrix([s,-2,(-1)/s],[0,1,1/s],[0,0,1/s]),
             /*matrix for adjoint action of y*/
             Y:matrix([s,0,0],[s*u,1,0],[(-s)*u^2,(-2)*u,1/s]),
             IX,/*inverse matrix of X*/
             IY,/*inverse matrix of Y*/
             S:ident(3),/*marix for series of W or W^{-1}*/
             AS:ident(3),/*adjoint matrix of S*/
             W,/*matrix W=[Y,IX] */
             IW,/*matrix [IX,Y]*/
             d:1,/*the highest degree of u in the numerator 
                    of R-torsion*/
             k:1,/*the highest degree of u in the Riley poly*/
             p:0,
             r:R(m),/*the Riley poly*/
             r1 /*a polynomial removed the top term of u 
                  from the Riley polynomial*/
            ],
            IX:invert(X),
            IY:invert(Y),
            W:Y.IX.IY.X,
            IW:IX.Y.X.IY,
            /*calculating the numerator of R-torsion*/
            if m>0 then 
              block(
                /*calculation of S*/
                for i:1 thru m-1 do(S:ident(3)+S.IW),
                AS:adjoint(S),
                /*the numerator of R-torsion*/
                p:p+mattrace(IX.S),
                p:p+3*mattrace(X.AS)+mattrace(IY.W.AS),
                p:p-mattrace(IX.S)*mattrace((ident(3)+IX.Y).S),
                p:p+mattrace(IX.S.(ident(3)+IX.Y).S),
                p:p+(2-s+(-1)/s)^2*determinant(S),
                /*The top term of u in the Riley poly r
                   is given by -u^(2m-1).
                  We use the relation u^(2m-1) = r + u^(2m-1) later*/
                k:2*m-1,/*the highets degree of u in the Riley poly*/
                r1:r+u^(2*m-1)
              )
            else 
              block(
                /*calculation of S*/
                for i:1 thru -m-1 do(S:ident(3)+S.W),
                AS:adjoint(S),
                p:p-mattrace(IX.S.W),
                p:p+3*mattrace(X.IW.AS)+mattrace(IY.AS),
                p:p-mattrace(IX.S.W)*mattrace((ident(3)+IX.Y).S.W),
                p:p+mattrace(IX.S.W.(ident(3)+IX.Y).S.W),
                p:p-(2-s+(-1)/s)^2*determinant(S),
                /*The top term of u in the Riley poly r 
                   is given by u^(2|m|). 
                  We use the relation u^(2|m|) = -r + u^(2|m|) later*/
                k:2*(-m),/*the highets degree of u in the Riley poly*/
                r1:-r+u^(2*(-m))
              ),
            p:expand(p),
            /* simplify by using r (decreasing the degrees of u)*/
            /* set the degree of u in p*/
            d:hipow(p,u),
            /*decreasing the degrees of u*/
            for j:1 while d >= k do(
              p:subst(r1*u^(d-k),u^d,p),
              p:expand(p),
              d:hipow(p,u)
            ),
            p:factor(p),
            /*multiplying p 
               by the denominator of twisted Alexander*/
            p:expand(p*(s/(s^2-2*s+1))),
            p:factorout(p,s),
            r:factorout(r,s),
            print("The Riley polynomial of J(2,",2*m,"):",r),
            print("The Reidemeister torsion for J(2,",2*m,"):"),
            return(p)
          )
      else print(m,"is not an integer.");
\end{verbatim}
\end{small}

\subsection{Tables}
%%%%%%%
%TABLES
%%%%%%%
In this appendix, except in the case of the trefoil knot (the only twist knot which is not hyperbolic), $u$ denotes the root of Riley's equation $\phi_{K}(1, u) = 0$ corresponding to the discrete and faithful representation of the complete hyperbolic structure.

Tables~\ref{Table+} and~\ref{Table-} give the non--abelian Reidemeister torsion for twist knots at the holonomy (except in the case of the trefoil) with respect to the corresponding root of Riley's equation and to the cusp shape of the knot exterior.

\hfill

%%
% Table of J(2, 2m) (1 \leq m \leq 10)
%%
\begin{table}[p]
%\centering 
%\begin{flushleft}
%\begin{footnotesize}
{\begin{tabular}[l]{|c|p{11.3cm}|}
\hline
\multicolumn{2}{|c|}
{\bf Non--abelian torsion for $J(2, 2m)$\, $(1\leq m\leq 10)$}\\
\hline
\Gcenter{3}{$m$} \rule{0in}{4mm}
   &Torsion at the holonomy ${\mathbb T}^{J(2, 2m)}_{\lambda}(\rho_u)$ (divided by a sign $-\tau_0$)\\[1mm]
\cline{2-2} \rule{0in}{4mm}
   &Result by substituting $u=4/(\cusp-2)$ into ${\mathbb T}^{J(2, 2m)}_{\lambda}(\rho_u)$, where $\cusp$ is the cusp shape\\[1mm]
\hline

\Gcenter{2}{$1$} \rule{0in}{3mm}
   &$3$\\
\cline{2-2} \rule{0in}{3mm}
   &$3$\\
\hline

\Gcenter{2}{$2$} \rule{0in}{3mm}
   &$13u^2-7u+19$\\
\cline{2-2} \rule{0in}{3mm}
   &$(19\cusp^2-104\cusp+340) / (\cusp-2)^2$\\
\hline

\Gcenter{2}{$3$} \rule{0in}{3mm}
   &$26u^4-17u^3+98u^2-45u+55$\\
\cline{2-2} \rule{0in}{3mm}
   &$(55\cusp^4-620\cusp^3+3968\cusp^2-11280\cusp+17424) / (\cusp-2)^4$\\
\hline

\Gcenter{2}{$4$} \rule{0in}{3mm}
   &$46u^6-34u^5+263u^4-157u^3+402u^2-159u+118$\\
\cline{2-2} \rule{0in}{3mm}
   &$2(59\cusp^6-1026\cusp^5+9936\cusp^4-52912\cusp^3+180592\cusp^2-352032\cusp+369280) / (\cusp-2)^6$\\
\hline

\Gcenter{3}{$5$} \rule{0in}{3mm}
   &$69u^8-54u^7+540u^6-366u^5+1360u^4-733u^3+1186u^2-411u+215$\\
\cline{2-2} \rule{0in}{3mm}
   &$(215\cusp^8-5084\cusp^7+66072\cusp^6-509040\cusp^5+2656960\cusp^4-9378624\cusp^3+22613632\cusp^2-33723648\cusp+26688768) / (\cusp-2)^8$\\
\hline

\Gcenter{3}{$6$} \rule{0in}{3mm}
   &$99u^{10}-81u^9+971u^8-710u^7+3400u^6-2123u^5+5052u^4-2469u^3+2875u^2-884u+353$\\
\cline{2-2} \rule{0in}{3mm}
   &$(353\cusp^{10}-10596\cusp^9+173188\cusp^8-1742080\cusp^7+12219808\cusp^6-61550208\cusp^5+228030592\cusp^4-612284416\cusp^3+1160955136\cusp^2-1411093504\cusp+903214080) / (\cusp-2)^{10}$\\
\hline

\Gcenter{5}{$7$} \rule{0in}{3mm}
   &$132u^{12}-111u^{11}+1566u^{10}-1203u^9+7057u^8-4810u^7+14996u^6-8647u^5+15044u^4-6710u^3+6076u^2-1678u+539$\\
\cline{2-2} \rule{0in}{3mm}
   &$(539\cusp^{12}-19648\cusp^{11}+387176\cusp^{10}-4799040\cusp^9+42208784\cusp^8-274741248\cusp^7+1363062528\cusp^6-5187840000\cusp^5+15118560512\cusp^4-33001455616\cusp^3+51856091136\cusp^2-53202206720\cusp+28544299008)/(\cusp-2)^{12}$\\
\hline

\Gcenter{5}{$8$} \rule{0in}{3mm}
   &$172u^{14}-148u^{13}+2383u^{12}-1899u^{11}+13098u^{10}-9475u^9+36258u^8-23106u^7+52884u^6-28275u^5+38518u^4-15774u^3+11636u^2-2914u+780$\\
\cline{2-2} \rule{0in}{3mm}
   &$4(195\cusp^{14}-8374\cusp^{13}+193288\cusp^{12}-2846448\cusp^{11}+30095568\cusp^{10}-239812000\cusp^9+1487434752\cusp^8-7287857664\cusp^7+28404952320\cusp^6-87772645888\cusp^5+212579837952\cusp^4-393068802048\cusp^3+529782681600\cusp^2-471281852416\cusp+218188021760) / (\cusp-2)^{14}$\\
\hline

\Gcenter{6}{$9$} \rule{0in}{3mm}
   &$215u^{16}-188u^{15}+3416u^{14}-2796u^{13}+22210u^{12}-16767u^{11}+76022u^{10}-51847u^9+146639u^8-87602u^7+157972u^6-78647u^5+87864u^4-33238u^3+20652u^2-4730u+1083$\\
\cline{2-2} \rule{0in}{3mm}
   &$(1083\cusp^{16}-53576\cusp^{15}+1417872\cusp^{14}-24177568\cusp^{13}+298484224\cusp^{12}-2810875520\cusp^{11}+20889506048\cusp^{10}-124832580096\cusp^9+606632721920\cusp^8-2406783375360\cusp^7+7784342106112\cusp^6-20354210914304\cusp^5+42341634637824\cusp^4-68068269064192\cusp^3+80435317243904\cusp^2-63142437978112\cusp+25688198283264) / (\cusp-2)^{16}$\\
\hline

\Gcenter{7}{$10$} \rule{0in}{3mm}
    &$265u^{18}-235u^{17}+4739u^{16}-3964u^{15}+35520u^{14}-27711u^{13}+144776u^{12}-103759u^{11}+348155u^{10}-224404u^9+501055u^8-281458u^7+417368u^6-194245u^5+183500u^4-64454u^3+34537u^2-7285u+1455$\\
\cline{2-2} \rule{0in}{3mm}
    &$(1455\cusp^{18}-81520\cusp^{17}+2433812\cusp^{16}-47158400\cusp^{15}+665730240\cusp^{14}-7230947840\cusp^{13}+62585931008\cusp^{12}-440831379456\cusp^{11}+2561522930176\cusp^{10}-12371911213056\cusp^9+49827680770048\cusp^8-167134091640832\cusp^7+464297966682112\cusp^6-1056316689612800\cusp^5+1931794260754432\cusp^4-2753296051208192\cusp^3+2901909811167232\cusp^2-2040504620417024\cusp+741196988416000)/(\cusp-2)^{18}$\\
\hline
\end{tabular}}
\caption{Table for the sequence of knots $J(2, 2m)$, $(1\leq m\leq 10)$.}
\label{Table+}
%\end{footnotesize}
%\end{flushleft}
\end{table}

%%
% Table of J(2, -2m) (1 \leq m \leq 10)
%%
\begin{table}[p]
%  \centering 
%\begin{flushleft}
%\begin{footnotesize}
{\begin{tabular}[l]{|c|p{11.3cm}|}
\hline
\multicolumn{2}{|c|}
{\bf Non--abelian torsion for $J(2, -2m)$\, $(1\leq m\leq 10)$}\\
\hline

\Gcenter{3}{$m$} \rule{0in}{4mm}
   &Torsion at the holonomy ${\mathbb T}^{J(2,-2m)}_{\lambda}(\rho_u)$ (divided by a sign $\tau_0$)\\[1mm]
\cline{2-2} \rule{0in}{4mm}
   &Result by substituting $u=4/(\cusp-2)$ into ${\mathbb T}^{J(2, -2m)}_{\lambda}(\rho_u)$, where $\cusp$ is the cusp shape \\[1mm]
\hline

\Gcenter{2}{$1$} \rule{0in}{3mm}
   &$-3$\\
\cline{2-2} \rule{0in}{3mm}
   &$-3$\\
\hline

\Gcenter{2}{$2$} \rule{0in}{3mm}
   &$7u^3+u^2+14u-5$\\
\cline{2-2} \rule{0in}{3mm}
   &$-(5\cusp^3-86\cusp^2+268\cusp-680)/(\cusp-2)^3$\\
\hline

\Gcenter{2}{$3$} \rule{0in}{3mm}
   &$17u^5+8u^4+79u^3+26u^2+73u+1$\\
\cline{2-2} \rule{0in}{3mm}
   &$(\cusp^5+282\cusp^4-1880\cusp^3+9488\cusp^2-22448\cusp+34848)/(\cusp-2)^5$\\
\hline

\Gcenter{2}{$4$} \rule{0in}{3mm}
   &$34u^7+22u^6+225u^5+119u^4+439u^3+162u^2+229u+22$\\
\cline{2-2}\rule{0in}{3mm}
   &$2(11\cusp^7+304\cusp^6-3276\cusp^5+25488\cusp^4-112432\cusp^3+359808\cusp^2-661824\cusp+738560)/(\cusp-2)^7$\\
\hline

\Gcenter{3}{$5$}\rule{0in}{3mm}
   &$54u^9+39u^8+474u^7+300u^6+1411u^5+730u^4+1619u^3+586u^2+551u+65$\\
\cline{2-2}\rule{0in}{3mm}
   &$(65\cusp^9+1034\cusp^8-16528\cusp^7+175520\cusp^6-1125280\cusp^5+5374144\cusp^4-17783040\cusp^3+42336768\cusp^2-62848768\cusp+53377536)/(\cusp-2)^9$\\
\hline

\Gcenter{4}{$6$} \rule{0in}{3mm}
   &$81u^{11}+63u^{10}+872u^9+611u^8+3462u^7+2104u^6+6167u^5+3036u^4+4714u^3+1615u^2+1129u+137$\\
\cline{2-2} \rule{0in}{3mm}
   &$(137\cusp^{11}+1502\cusp^{10}-34340\cusp^9+468616\cusp^8-3940960\cusp^7+25007808\cusp^6-117308544\cusp^5+420442368\cusp^4-1109472000\cusp^3+2123257344\cusp^2-2628703232\cusp+1806428160)/(\cusp-2)^{11}$\\
\hline
 
\Gcenter{5}{$7$} \rule{0in}{3mm}
   &$111u^{13}+90u^{12}+1425u^{11}+1062u^{10}+7105u^9+4747u^8+17286u^7+9956u^6+21045u^5+9752u^4+11610u^3+3724u^2+2070u+245$\\
\cline{2-2} \rule{0in}{3mm}
   &$(245\cusp^{13}+1910\cusp^{12}-62696\cusp^{11}+1057552\cusp^{10}-11025808\cusp^9+87196704\cusp^8-526180096\cusp^7+2499806720\cusp^6-9272200448\cusp^5+26825366016\cusp^4-58773907456\cusp^3+94191718400\cusp^2-99494326272\cusp+57088598016)/(\cusp-2)^{13}$\\
\hline

\Gcenter{5}{$8$} \rule{0in}{3mm}
   &$148u^{15}+124u^{14}+2195u^{13}+1711u^{12}+13125u^{11}+9354u^{10}+40453u^9+25698u^8+68070u^7+37044u^6+60745u^5+26422u^4+25394u^3+7604u^2+3502u+396$\\
\cline{2-2} \rule{0in}{3mm}
   &$4(99\cusp^{15}+532\cusp^{14}-26060\cusp^{13}+529856\cusp^{12}-6606160\cusp^{11}+62595520\cusp^{10}-460836032\cusp^9+2719805952\cusp^8-12910096128\cusp^7+49528327168\cusp^6-152285103104\cusp^5+370984124416\cusp^4-695663718400\cusp^3+961422573568\cusp^2-884915453952\cusp+436376043520)/(\cusp-2)^{15}$\\
\hline

\Gcenter{8}{$9$} \rule{0in}{3mm}
   &$188u^{17}+161u^{16}+3172u^{15}+2552u^{14}+22171u^{13}+16540u^{12}+82961u^{11}+56366u^{10}+179181u^9+108029u^8+224190u^7+115204u^6+154037u^5+62916u^4+50650u^3+14172u^2+5570u+597$\\
\cline{2-2} \rule{0in}{3mm}
   &$(597\cusp^{17}+1982\cusp^{16}-161440\cusp^{15}+3885760\cusp^{14}-56503104\cusp^{13}+624108160\cusp^{12}-5417133568\cusp^{11}+38142084096\cusp^{10}-219869856256\cusp^9+1045959103488\cusp^8-4105465389056\cusp^7+13257704030208\cusp^6-34864687169536\cusp^5+73473552449536\cusp^4-120432383098880\cusp^3+146315203575808\cusp^2-119078046597120\cusp+51376396566528)/(\cusp-2)^{17}$\\
\hline

\Gcenter{8}{$10$} \rule{0in}{3mm}
    &$235u^{19}+205u^{18}+4434u^{17}+3659u^{16}+35404u^{15}+27360u^{14}+155687u^{13}+111176u^{12}+410964u^{11}+266255u^{10}+665399u^9+380935u^8+647346u^7+314408u^6+353985u^5+135980u^4+94041u^3+24637u^2+8440u+855$\\
\cline{2-2} \rule{0in}{3mm}
    &$(855\cusp^{19}+1270\cusp^{18}-236348\cusp^{17}+6649256\cusp^{16}-110706240\cusp^{15}+1397436800\cusp^{14}-13965196032\cusp^{13}+114155938304\cusp^{12}-773068443136\cusp^{11}+4380077954048\cusp^{10}-20842143467520\cusp^9+83401615028224\cusp^8-279834167951360\cusp^7+782314999349248\cusp^6-1800640172326912\cusp^5+3349441682210816\cusp^4-4881403696185344\cusp^3+5296224268058624\cusp^2-3862939033141248\cusp+1482393976832000)/(\cusp-2)^{19}$\\
\hline
\end{tabular}}
\caption{Table for the sequence of knots $J(2, -2m)$, $(1\leq m\leq 10)$.}
\label{Table-}
%\end{footnotesize}
%\end{flushleft}
\end{table}

\clearpage

%%%%%%%%%%%%%%%%%%%%
% reference 
%%%%%%%%%%%%%%%%%%%%
\section*{Acknowledgments}

The first author  is supported by the {European Community} with Marie Curie Intra--European Fellowship (MEIF--CT--2006--025316). While writing the paper, J.D. visited the CRM. He thanks the CRM for hospitality.
The third author is partially supported by 
the 21st century COE program at Graduate School of Mathematical Sciences, 
University of Tokyo.
On finishing the paper, J.D. visited the Department of Mathematics of Tokyo Institute of Technology.
J.D. and Y.Y. wish to thank Hitoshi Murakami and
 TiTech's Department of Mathematics for invitation and hospitality. The authors also want to thank Joan Porti for his comments and remarks.

\bibliographystyle{amsalpha}
\bibliography{twist_knots}

\newcommand{\noop}[1]{}
\providecommand{\bysame}{\leavevmode\hbox to3em{\hrulefill}\thinspace}
\providecommand{\MR}{\relax\ifhmode\unskip\space\fi MR }
\providecommand{\MRhref}[2]{%
  \href{http://www.ams.org/mathscinet-getitem?mr=#1}{#2}
}
\providecommand{\href}[2]{#2}
\begin{thebibliography}{CGLS87}

\bibitem[Ada94]{Adams}
C.~Adams, \emph{The knot book}, W. H. Freeman and Company, 1994.

\bibitem[Cal06]{Calegari}
D.~Calegari, \emph{Real places and torus bundles}, Geom. Dedicata
  \textbf{\textbf{118}} (2006), 209--227.

\bibitem[CGLS87]{DehnSurgery}
M.~Culler\noop{1}, C.~McA. Gordon, J.~Luecke, and P.~Shalen, \emph{Dehn surgery
  on knots}, Ann. of Math. \textbf{\textbf{125}} (1987), 237--300.

\bibitem[CS83]{CS:1983}
M.~Culler\noop{2} and P.~Shalen, \emph{Varieties of group representations and
  splittings of $3$-manifolds}, Ann. of Math. \textbf{\textbf{117}} (1983),
  109--146.

\bibitem[DK07]{DK05}
J.~Dubois and R.~Kashaev, \emph{On the asymptotic expansion of the colored
  {J}ones polynomial for torus knots}, Math. Ann. \textbf{339} (2007),
  757--782.

\bibitem[Dub05]{JDFourier}
J.~Dubois, \emph{Non abelian {R}eidemeister torsion and volume form on the
  {${\mathrm{SU}}(2)$}-representation space of knot groups}, Ann. Institut
  Fourier \textbf{\textbf{55}} (2005), 1685--1734.

\bibitem[Dub06]{JDFibre}
\bysame, \emph{Non abelian twisted {R}eidemeister torsion for fibered knots},
  Canad. Math. Bull. \textbf{\textbf{49}} (2006), 55--71.

\bibitem[FV07]{FV07}
S.~Friedl and S.~Vidussi, \emph{Symplectic $\text{S}^1 \times {N}^3$, subgroup
  separability, and vanishing {T}hurston norm}, 2007, preprint
  arXiv:math/0701717.

\bibitem[GKM05]{GKM}
H.~Goda, T.~Kitano, and T.~Morifuji, \emph{{R}eidemeister torsion, twisted
  {A}lexander polynomial and fibered knots}, {Comment. Math. Helv.} \textbf{80}
  (2005), 51--61.

\bibitem[Heu03]{Heu:2003}
M.~Heusener, \emph{{An orientation for the $\mathrm{SU}(2)$-representation
  space of knot groups}}, Topology and its Applications \textbf{127} (2003),
  175--197.

\bibitem[HS04]{HS}
J.~Hoste and P.~D. Shanahan, \emph{A formula for the {A}-polynomial of twist
  knots}, J. Knot Theory Ramifications \textbf{\textbf{13}} (2004), 193--209.

\bibitem[Kit96]{Kitano}
T.~Kitano, \emph{Twisted {A}lexander polynomial and {R}eidemeister torsion},
  Pacific J. Math. \textbf{\textbf{174}} (1996), 431--442.

\bibitem[KL99]{KL}
P.~Kirk and C.~Livingston, \emph{Twisted {A}lexander {I}nvariants,
  {R}eidemeister torsion, and {C}asson-{G}ordon invariants}, Topology
  \textbf{\textbf{38}} (1999), 635--661.

\bibitem[Le94]{Le}
T.~Le, \emph{Varieties of representations and their subvarieties of cohomology
  jumps for certain knot groups}, Russian Acad. Sci. Sb. Math.
  \textbf{\textbf{78}} (1994), 187--209.

\bibitem[Lei02]{Leininger}
C.~J. Leininger, \emph{Surgeries on one component of the {W}hitehead link are
  virtually fibered}, Topology \textbf{\textbf{41}} (2002), 307--320.

\bibitem[Lin01]{Lin01}
X-S. Lin, \emph{Representations of knot groups and twisted {A}lexander
  polynomials}, Acta Math. Sin. (Engl. Ser.) \textbf{17} (2001), no.~3,
  361--380.

\bibitem[Men84]{Menasco}
W.~Menasco, \emph{Closed incompressible surfaces in alternating knot and link
  complements}, Topology \textbf{\textbf{23}} (1984), 37--44.

\bibitem[Mer03]{Merris:2003}
R.~Merris, \emph{Combinatorics}, 2nd ed., John Wiley and Sons, 2003.

\bibitem[Mil66]{Milnor:1966}
J.~Milnor, \emph{Whitehead torsion}, Bull. Amer. Math. Soc.
  \textbf{\textbf{72}} (1966), 358--426.

\bibitem[Por97]{Porti:1997}
J.~Porti, \emph{Torsion de {R}eidemeister pour les vari\'et\'es hyperboliques},
  vol. 128, Mem. Amer. Math. Soc., no. 612, AMS, 1997.

\bibitem[Ril84]{Riley}
R.~Riley, \emph{Nonabelian representations of $2$-bridge knot groups}, Quart.
  J. Math. Oxford Ser. $(2)$ \textbf{\textbf{35}} (1984), 191--208.

\bibitem[Rol90]{Rolfsen}
D.~Rolfsen, \emph{Knots and links}, Publish or Perish Inc., 1990.

\bibitem[Sa06]{Max}
W.~Schelter and al, \emph{Maxima}, 2006, available at
  \texttt{http://maxima.sourceforge.net/}.

\bibitem[Sha02]{HandBookGT}
P.~B. Shalen, \emph{Representations of $3$-manifold groups}, Handbook of
  Geometric Topology (R.~J. Daverman and R.~B. Sher, eds.), North-Holland
  Amsterdam, 2002, pp.~955--1044.

\bibitem[SW06]{SW06}
D.~Silver and S.~Williams, \emph{Twisted {A}lexander polynomials detect the
  unknot}, Algebr. Geom. Topol. \textbf{\textbf{6}} (2006), 1903--1923.

\bibitem[Thu02]{TN}
W.~Thurston, \emph{The {G}eometry and {T}opology of {T}hree-{M}anifolds},
  Electronic Notes available at
  \texttt{http://www.msri.org/publications/books/gt3m/}, 2002.

\bibitem[Tur01]{Turaev:2000}
V.~Turaev, \emph{Introduction to combinatorial torsions}, Lectures in
  Mathematics, Birkh{\"a}user, 2001.

\bibitem[Tur02]{Turaev:2002}
\bysame, \emph{Torsions of $3$-dimensional manifolds}, Progress in Mathematics,
  vol. \textbf{208}, Birkh{\"a}user, 2002.

\bibitem[Wad94]{Wada94}
M.~Wada, \emph{Twisted {A}lexander polynomial for finitely presentable groups},
  Topology \textbf{33} (1994), no.~2, 241--256.

\bibitem[Wal78]{W}
F.~Waldhausen, \emph{Algebraic {K}-theory of generalized free products {I},
  {II}.}, Ann. of Math. \textbf{\textbf{108}} (1978), 135--204.

\bibitem[Wee99]{snappea}
J.~Weeks, \emph{Snap{P}ea}, 1999, available at
  \texttt{http://www.geometrygames.org/SnapPea/}.

\bibitem[Yam05]{YY1}
Y.~Yamaguchi\noop{1}, \emph{A relationship between the non-acyclic
  {R}eidemeister torsion and a zero of the acyclic {R}eidemeister torsion},
  2005, to appear in {\it Ann. Institut Fourier} (arXiv:math.GT/0512267).

\bibitem[Yam07]{YY2}
Y.~Yamaguchi\noop{2}, \emph{Limit values of the non-acyclic {R}eidemeister
  torsion for knots}, Algebr. Geom. Topol. \textbf{\textbf{7}} (2007),
  1485--1507.

\end{thebibliography}

\end{document}